\documentclass[11pt,reqno]{amsart}
\usepackage{multicol}
\usepackage{indentfirst, amssymb,amsmath,amsthm, mathrsfs,cite}
\usepackage{xcolor}
\newtheorem*{theoA}{Theorem A}
\newtheorem*{theoB}{Theorem B}

\newtheorem{theorem}{Theorem}[section]
\newtheorem{lemma}{Lemma}[section]

\newtheorem{example}{Example}[section]

\newtheorem*{Canonical function}{Canonical function}
\newcommand{\pa}{\partial}
\newcommand{\ol}{\overline}
\newcommand{\be}{\begin{equation}}
\newcommand{\ee}{\end{equation}}
\newcommand{\beas}{\begin{eqnarray*}}
\newcommand{\eeas}{\end{eqnarray*}}
\newcommand{\bea}{\begin{eqnarray}}
\newcommand{\eea}{\end{eqnarray}}
\newcommand{\bs}{\begin{small}}
\newcommand{\es}{\end{small}}
\renewcommand{\epsilon}{\varepsilon}
\numberwithin{equation}{section}

\begin{document}
\title[Entire functions of several]{Entire functions of several complex variables satisfying certain Fermat-type PDDE\lowercase{s}}
\author[H. Y. XU, R. Mandal and R. Biswas]{Hong Yan Xu, Rajib Mandal and Raju Biswas}
\date{}
\address{School of Arts and Sciences, Suqian University, Suqian, Jiangsu 223800, China.}
\email{xhyhhh@126.com}
\address{Department of Mathematics, Raiganj University, Raiganj, West Bengal-733134, India.}
\email{rajibmathresearch@gmail.com}
\address{Department of Mathematics, Raiganj University, Raiganj, West Bengal-733134, India.}
\email{rajubiswasjanu02@gmail.com}
\maketitle
\let\thefootnote\relax
\footnotetext{2020 Mathematics Subject Classification: 39A45, 39A14, 39B32, 32W50, 30D35.}
\footnotetext{Key words and phrases: Fermat-type equation, Several complex variables, Partial differential-difference equation, Nevanlinna theory.}
\footnotetext{Type set by \AmS -\LaTeX}
\begin{abstract} In this paper, we solve certain Fermat-type partial differential-difference equations for finite order entire functions of several complex variables. These 
results are significant generalizations of some earlier findings, especially those of Haldar and Ahamed (Entire solutions of several quadratic binomial and trinomial partial differential-difference equations in $\mathbb{C}^2$, Anal. Math. Phys., 12 (2022)).
In addition, the results improve the previous results from the situation with two complex variables to the situation with several complex variables. To support our results, we have included several examples.
\end{abstract}
\section{Introduction, Definitions and Results}
A meromorphic function on the $n-$dimensional complex space is a function that can be expressed as a quotient of two holomorphic functions on the same space, without any zero sets in common. Notationally, we write $f:=\frac{g}{h}$, where $g$ and $h$ are relatively prime holomorphic functions on $\mathbb{C}^n$ such that $h \not\equiv 0$ and 
$f^{-1}(\infty)\not=\mathbb{C}^n$. In particular, the entire function of several complex variables are holomorphic throughout $\mathbb{C}^n$. \\[2mm]
\indent Let $z=(z_1,z_2,\ldots,z_n)\in\mathbb{C}^n$, $a\in\mathbb{C}\cup\{\infty\}$, $k\in\mathbb{N}$ and $r>0$. We have adapted some of the notations from \cite{14,S1,Y1}. Let $\ol B_n(r):=\{z\in\mathbb{C}^n: |z|\leq r\}$, where $|z|^2:=\sum_{j=1}^n|z_j|^2$. 
The exterior derivative splits $d:=\pa+\ol{\pa}$ and twists to $d^c:=\frac{i}{4\pi}(\ol\pa-\pa)$. The standard Kaehler metric on $\mathbb{C}^n$ is 
given by $v_n(z):=dd^c|z|^2$. Define $\omega_n(z):=dd^c\log |z|^2\geq 0$ and $\sigma_n(z):=d^c \log |z|^2\wedge \omega_n^{n-1}(z)$ on $\mathbb{C}^n\setminus\{0\}$. Thus $\sigma_n(z)$ defines a positive measure on $\pa B_n:=\{z\in\mathbb{C}^n: |z|= r\}$ with total measure $1$.
The zero-multiplicity of a holomorphic function $h$ at a point 
$z\in \mathbb{C}^n$ is defined to be the order of vanishing of $h$ at $z$ and denoted by $\mathcal{D}_h^0(z)$. A divisor of $f$ on $\mathbb{C}^n$ is an integer valued function which is locally the difference between the 
zero-multiplicity functions of $g$ and $h$ and it is denoted by $\mathcal{D}_f:=\mathcal{D}_g^0-\mathcal{D}_h^0$ (see P. 381, \cite{51}). Let $a\in\mathbb{C}\cup\{\infty\}$ be such that $f^{-1}(a)\not=\mathbb{C}^n$. Then the $a$-divisor $\nu_f^a$ of $f$ 
is the divisor associated with the holomorphic functions $g-ah$ and $h$ (see P. 346, \cite{14} and P. 12, \cite{12}). 
In \cite{Y1}, Ye has defined the counting function and the valence function with respect to $a$ respectively as follows:
\beas n(r,a,f):=r^{2-2n}\int_{S(r)} \nu_f^a v_n^{n-1}\;\text{and}\;N(r,a,f):=\int_0^r\frac{ n(r,a,f)}{t}dt. \eeas
We write 
\bs\beas N(r,a,f)=\left\{\begin{array}{llc}
&N\left(r,\frac{1}{f-a}\right),&\text{when}\;a\not=\infty\\[2mm]
&N(r,f),&\text{when}\; a=\infty.\end{array}\right.\eeas \es
The proximity function \cite{14,Y1} of $f$ is defined as follows :
\beas\left\{\begin{array}{lll}
&m(r,f):=\int_{\pa B_n(r)} \log^+|f(z)|\sigma_n(z),&\text{when}\;a=\infty\\[2mm]
&m\left(r, \frac{1}{f-a}\right):=\int_{\pa B_n(r)} \log^+\frac{1}{|f(z)-a|}\sigma_n(z),&\text{when}\;a\not=\infty.\end{array}\right.\eeas
By denoting $S(r):= \ol{B}_n(r)\cap \;\text{supp}\;\nu_f^a$, where 
$\text{supp}\;\nu_f^a=\ol{\left\{z\in\mathbb{C}^n: \nu_f^a(z)\not=0\right\}}$ (see P. 346, \cite{14}). The notation $N_k\left(r,\frac{1}{f-a}\right)$ is known as truncated valence function. In particular, $N_1\left(r,\frac{1}{f-a}\right)=\ol N\left(r,\frac{1}{f-a}\right)$ is the truncated valence function of simple $a$-divisors of $f$ in $S(r)$. 
In $N_k\left(r,\frac{1}{f-a}\right)$, the $a$-divisors of $f$ in $S(r)$ of multiplicity $m$  are counted $m$-times if $m< k$ and $k$-times if $m\geq k$. The Nevanlinna characteristic function is defined by $T(r,f)=N(r,f)+m(r,f)$, which is increasing for $r$. The order of a meromorphic function $f$ is denoted by $\rho(f)$ and is defined by 
\beas \rho(f)=\varlimsup\limits_{r\to \infty}\frac{\log^+ T(r,f)}{\log r},\;\;\text{where}\;\log^{+}x=\max\{\log x,0\}.\eeas
Recall that $T(r,\alpha)=S(r,f)$ indicates that the meromorphic function $\alpha$ is a small function of $f$, where $S(r,f)$ is any quantity that satisfies 
$S(r,f)=o(T(r,f))$ as $r\rightarrow \infty$, possibly outside of an exceptional set $E$ of finite linear measure $(\int_E dr<+\infty)$. 
For additional information, see \cite{51,12i,12,15,19,21,S1,Y1} and its references.\\[2mm]
\indent We now consider the Fermat-type equation
\bea\label{eq1.1} f^n(z)+g^n(z)=1,\;\text{where}\;n\in\mathbb{N}.\eea 
For solutions of the equation (\ref{eq1.1}) on $\mathbb{C}$, the classical results are found in \cite{7,49,39,17,50,7}.
During the last two decades, an increasing number of researchers have expressed interest in studying the Fermat-type equations for entire and meromorphic solutions by taking some 
variation of (\ref{eq1.1}). Yang and Li \cite{40} were the pioneers for introducing the study on transcendental meromorphic solutions of Fermat-type differential equations on 
$\mathbb{C}$. The first researcher who investigated meromorphic solutions of the Fermat-type difference equations and Fermat-type differential-difference equations on 
$\mathbb{C}$  was Liu \cite{Liu1}.  The classical results for solutions of the equation (\ref{eq1.1}) on several complex variables are found in  \cite[\textrm{Theorem 2.3}]{700} and \cite[\textrm{Theorem 1.3}]{20}. For the most recent and leading improvements in these directions, we refer to \cite{Ban1,601,7,49,28,16,Liu2,602,603,604,600,Raj2} and the references therein.\\[2mm]
\indent An equation that contains partial derivatives of $f$ is referred to as a partial differential equation (in brief, PDE). If the equation also includes shifts or differences of $f$, it is referred to as a partial differential-difference equation (in brief, PDDE).\\[2mm]
\indent In recent years, researchers have focused their attention on the investigation of solutions to Fermat-type PDDEs. Let 
\bea\label{esl}\sum_{i=1}^n\left(\frac{\pa u}{\pa z_i}\right)^m=1\eea
be the certain non-linear first order PDE introducing from the analogy with the Fermat-type equation $\sum_{i=1}^n\left(f_i\right)^m=1$, where $u:\mathbb{C}^n\to \mathbb{C}$, $z_i\in\mathbb{C}$, $f_i:\mathbb{C}\to \mathbb{C}$, and $m,n\geq 2$.
In 1999, Saleeby \cite{27} first started to study about the solutions of the Fermat-type PDEs and obtained the results for entire solutions of (\ref{esl}) on $\mathbb{C}^2$. 
Afterwards, in 2004, Li \cite{28} extended these results to $\mathbb{C}^n$. In 2008, Li \cite{16} considered the equation (\ref{eq1.1}) with $n=2$ and showed that meromorphic solutions $f$ and $g$ of that equation on 
$\mathbb{C}^2$ must be constant if and only if $\frac{\pa}{\pa z_2} f(z_1,z_2)$ and $\frac{\pa}{\pa z_1} g(z_1,z_2)$ have the same zeros (counting multiplicities). 
If $f=\frac{\pa}{\pa z_1} u$ and $g=\frac{\pa }{\pa z_2}u$, then any entire solutions of the PDE $f^2+g^2=1$ on $\mathbb{C}^2$ are necessarily linear \cite{13}.\\[2mm]
\indent The following results were established by Xu and Cao \cite{23,26} in 2018, who were the first to take consideration of both differential and difference operators in the Fermat-type equations on $\mathbb{C}^2$.
\begin{theoA} Let $c=(c_1,c_2)\in\mathbb{C}^2$. Then any transcendental entire solution with finite order of the Fermat-type PDDE 
$\left(\frac{\pa f(z)}{\pa z_1}\right)^2+f^2(z+c)=1$ has the form of $f(z_1,z_2)=\sin (Az_1+Bz_2+H(z_2))$, where $A,B$ are constants on $\mathbb{C}$ satisfying $A^2=1$ 
and $Ae^{i(Ac_1+Bc_2)}=1$ and $H(z_2)$ is a polynomial in one variable $z_2$ such that $H(z_2)\equiv H(z_2+c_2)$. In the special case whenever $c_2\not=0$, we have $f(z_1,z_2)=\sin(Az_1+Bz_2+\text{constant})$.\end{theoA}
\begin{theoB} Any transcendental entire solution with finite order of the Fermat type PDE \bea\label{esp} f^2(z)+\left(\frac{\pa f(z_1,z_2)}{\pa z_1}\right)^2=1\eea has 
the form of $f(z_1,z_2)=\sin(z_1+h(z_2))$, where $h(z_2)$ is a polynomial in $z_2$.\end{theoB}
By taking into consideration certain variations of the PDE (\ref{esp}), Chen and Xu \cite{35} improved upon the results of Xu and Cao \cite{23, 26} in 2021. The authors \cite{35} considered the 
following Fermat-type PDE
\bs\beas&&\left(a_2\frac{\pa f(z)}{\pa z_1}\right)^2+\left(a_3f(z)+a_4\frac{\pa^2f(z)}{\pa z_1^2}\right)^2=1,\\
&&\left(a_1 f(z)+a_2\frac{\pa f(z)}{\pa z_1}\right)^2+\left(a_3 f(z)+a_4\frac{\pa f(z)}{\pa z_2}\right)^2=1,\\
&&f^2(z)+\left(b_1\frac{\pa f(z)}{\pa z_1}+b_2\frac{\pa^2 f(z)}{\pa z_1^2}\right)^2=1\\\text{and}
&&f^2(z)+\left(b_1\frac{\pa f(z)}{\pa z_1}+b_2\frac{\pa^2 f(z)}{\pa z_1\pa z_2}\right)^2=1\eeas\es
and obtained the existence and forms of solutions for the transcendental entire functions of finite order of two complex variables $z_1$ and $z_2$.\\[2mm]
\indent In this direction, several results are found about the existence and forms of entire solutions of the Fermat-type PDEs and PDDEs on $\mathbb{C}^2$ (see \cite{41,500}).
In 2022, Haldar and Ahamed \cite{42} implemented many variations of the PDEs examined by Chen and Xu \cite{35}, substituting $f(z)$ with the difference function $f(z+c)$ and 
the difference operator $\Delta f(z)$ with constant coefficients. The following Fermat-type PDDEs were considered by the authors \cite{42}:
\bs\beas&&\left(a_2\frac{\pa f(z)}{\pa z_1}\right)^2+\left(a_3f(z+c)+a_4\frac{\pa^2f(z)}{\pa z_1^2}\right)^2=1,\\
&&\left(a_2\frac{\pa f(z)}{\pa z_1}\right)^2+\left(a_3\Delta f(z)+a_4\frac{\pa^2f(z)}{\pa z_1^2}\right)^2=1,\\
&&\left(a_1\Delta f(z)+a_2\frac{\pa f(z)}{\pa z_1}\right)^2+\left(a_3 \Delta f(z)+a_4\frac{\pa f(z)}{\pa z_1}\right)^2=1\\
&&f^2(z+c)+\left(b_1\frac{\pa f(z)}{\pa z_1}+b_2\frac{\pa^2 f(z)}{\pa z_1^2}\right)^2=1\\
\text{and}&&f^2(z+c)+\left(b_1\frac{\pa f(z)}{\pa z_1}+b_2\frac{\pa^2 f(z)}{\pa z_1\pa z_2}\right)^2=1\eeas\es
and obtained the forms of the solutions for the finite order transcendental entire functions.\\[2mm]
We refer to \cite{MB1,Raj10,35,Hal1,Xu2,Raj1,22,500,Xu1} and the references therein for the most recent developments in the aforementioned directions.
\section{The main results}
Let $I=(i_1,i_2,\ldots,i_n)\in\mathbb{Z}^n_+$ be a multi-index with length $\Vert I\Vert=\sum_{j=1}^n i_j$. 
We can write $\mathcal{P}(z)=\sum_{\Vert I\Vert=0}^d a_I z_1^{i_1}\cdots z_n^{i_n}$ to represent any polynomial of degree $d$ on $\mathbb{C}^n$, where $a_{I}\in\mathbb{C}$ so that $a_{I}$ are not all zero simultaneously for $\Vert I\Vert=d$.
\indent Let $c\in\mathbb{C}^n\setminus\{(0,0,\ldots,0)\}$. Note that for any monomial $H(z)$ of several complex variables with $\deg(H(z))\geq 2$, $H(z+c)-H(z)\not \equiv \text{constant}$. But $H(z+c)- H(z)\equiv$ constant is possible when the terms in $H(z)$ have powers of linear combinations of the variables with certain restrictions.\\[2mm] 
Suppose that for any $c\in\mathbb{C}^n\setminus\{(0,0,\ldots,0)\}$, $\mathcal{P}(z+c)-\mathcal{P}(z)\equiv B$, where $B\in\mathbb{C}$. Let $\mathcal{P}(z)=\sum_{j=1}^n a_jz_j+G(z)+A$, where $A\in\mathbb{C}$, $\deg (G(z))\geq 2$. Now, $\mathcal{P}(z+c)-\mathcal{P}(z)\equiv B$ implies that $\sum_{j=1}^n a_j c_j+G(z+c)-G(z)\equiv B$. Thus, we have $G(z+c)\equiv G(z)$ and $\sum_{j=1}^n a_jc_j=B$.\\[2mm]
\indent If $\deg(G(z))\geq 2$, then we can rewrite the polynomial $G(z)$ in such a way that it contains terms from the polynomials like 
$\Phi(t_{j_1}z_{j_1}+t_{j_2}z_{j_2}+\ldots+t_{j_m} z_{j_m})$ of $t_{j_1}z_{j_1}+t_{j_2}z_{j_2}+\ldots+t_{j_m} z_{j_m}$ such that $t_{j_1}c_{j_1}+t_{j_2}c_{j_2}+\ldots+t_{j_m} c_{j_m}=0$, $t_{j_1}$, $\cdots$, $t_{j_m}\in\mathbb{C}$ ($1\leq j_1,j_2,\cdots,j_m\leq n$) and $\deg(\Phi)\geq 2$.  Then, it is easy to se that $G(z)$ is periodic. In general, we can express $G(z)$ as 
\bea\label{K1} G(z)=\sum_{\lambda}G_\lambda(z)\quad\text{and}\quad G_\lambda(z)=\prod_{\alpha} G_\alpha (z),\eea
where $\lambda$ belongs to the finite index set $I_1$ of the family $\{G_\lambda(z) : \lambda \in I_1\}$ and $\alpha$ belongs to the finite index set $I_2$ of the family 
$\{G_\alpha(z) : \alpha \in I_2\}$ with
\beas G_\alpha(z)&=&\sum_{\substack{j_1,j_2=1,\\j_1<j_2}}^n \Phi_{2,\alpha,j_1,j_2}(\eta_{j_1}z_{j_1}+\eta_{j_2}z_{j_2})+\sum_{\substack{j_1,j_2,j_3=1,\\j_1<j_2<j_3}}^n \Phi_{3,\alpha,j_1,j_2,j_3}(\zeta_{j_1}z_{j_1}+\zeta_{j_2}z_{j_2}+\zeta_{j_3} z_{j_3})\nonumber\\
&&+\cdots+\sum_{\substack{j_1,j_2,\ldots,j_n=1,\\j_1<j_2<\ldots<j_n}}^n \Phi_{n,\alpha,j_1,j_2,\ldots,j_n}(t_{j_1}z_{j_1}+t_{j_2}z_{j_2}+\cdots+t_{j_n} z_{j_n})\nonumber\eeas
where $\eta_i,\zeta_i,t_i, A\in\mathbb{C}$ $(1\leq i\leq n)$, $\deg G(z)=\deg \mathcal{P}(z)$, $\Phi_{m,\alpha,j_1,j_2,\ldots,j_m}(t_{j_1}z_{j_1}+t_{j_2}z_{j_2}+\ldots+t_{j_m} z_{j_m})$ is a polynomial in $t_{j_1}z_{j_1}+t_{j_2}z_{j_2}+\ldots+t_{j_m} z_{j_m}$. Here $\eta_i,\zeta_i,t_i\in\mathbb{C}$ $(1\leq i\leq n)$ are chosen from the conditions $\eta_{j_1}c_{j_1}+\eta_{j_2}c_{j_2}=0$, $\zeta_{j_1}c_{j_1}+\zeta_{j_2}c_{j_2}+\zeta_{j_3} c_{j_3}=0$, $t_{j_1}c_{j_1}+t_{j_2}c_{j_2}+\ldots+t_{j_m} c_{j_m}=0$. It is also clear that for $ j_1=1,2,\ldots,n$, we have
\beas\frac{\pa G_\alpha(z)}{\pa z_{j_1}}=\eta_{j_1}\sum_{\substack{j_1,j_2=1,\\j_1<j_2}}^n \Phi_{2,\alpha,j_1,j_2}'(\eta_{j_1}z_{j_1}+\eta_{j_2}z_{j_2})+\zeta_{j_1}\sum_{\substack{j_1,j_2,j_3=1,\\j_1<j_2<j_3}}^n \Phi_{3,\alpha,j_1,j_2,j_3}'(\zeta_{j_1}z_{j_1}+\zeta_{j_2}z_{j_2}+\zeta_{j_3} z_{j_3})\\
+\ldots+t_{j_1}\sum_{\substack{j_1,j_2,\ldots,j_n=1,\\j_1<j_2<\ldots<j_n}}^n \Phi_{n,\alpha,j_1,j_2,\ldots,j_n}'(t_{j_1}z_{j_1}+t_{j_2}z_{j_2}+\cdots+t_{j_n} z_{j_n}).\hspace{5cm}\eeas
\indent The results of \cite{35,42,41,23,26,500} encourage us to go further on the generalization and improvement of the results. There are different ways to do this.\\
(I) Substituting the non-zero term "$e^{g(z)}$" for "$1$" on the R.H.S. of the equations, where $g(z)$ is a non-constant polynomial;\\
(II) Considering more general forms of partial differential-difference parts on the L.H.S. of the equations;\\
(III) Consider the aforementioned facts (I) - (II) combined on $\mathbb{C}^n$.\\[2mm]
\indent Note that, one can replace $1$ on the R.H.S. in the PDEs and PDDEs in \cite{35,42,41,23,26,500} by any finite order entire function, like, $F(z)$ such that $F(0)\not=0$ and attempt 
to generalize as well as improve the results in \cite{35,42,41,23,26,500}. Actually, for this kind of functions, we can easily apply \textrm{Lemma \ref{lem2}}, given in the lemma section. \\[2mm]
\indent In this research, keeping all the preceding facts in mind, we are investigating the existence of solutions of the following Fermat-type PDDEs on $\mathbb{C}^n$ for $1\leq \mu<\nu\leq n$:
\bea\label{e1}&&\left(a_1\frac{\pa f(z)}{\pa z_\mu}\right)^2+\left(a_2f(z)+a_3f(z+c)+a_4\frac{\pa^2 f(z)}{\pa z_\mu^2}\right)^2=e^{g(z)},\\[2mm]
\label{e2}&&\left(a_1\Delta f(z)+a_2\frac{\pa f(z)}{\pa z_\mu}\right)^2+\left(a_3\Delta f(z)+a_4\frac{\pa f(z)}{\pa z_\nu}\right)^2=e^{g(z)},\\[2mm]
\label{e3}&&a_1^2f^2(z+c)+\left(a_2\frac{\pa f(z)}{\pa z_\mu}+a_3\frac{\pa^2 f(z)}{\pa z_\mu^2}\right)^2=e^{g(z)}\\[2mm]\text{and}
\label{e4}&&a_1^2f^2(z+c)+\left(a_2\frac{\pa f(z)}{\pa z_\mu}+a_3\frac{\pa^2 f(z)}{\pa z_\mu\pa z_\nu}\right)^2=e^{g(z)},\eea
where $f(z)$ is a finite order transcendental entire function of several complex variables, $g(z)$ is a non-constant polynomial on $\mathbb{C}^n$ and $a_1,a_2,a_3,a_4\in\mathbb{C}\setminus\{0\}$. Throughout this paper, we denote
\beas\begin{array}{lll}
 y =\left(z_1,z_2,\dots, z_{\mu-1},a_1a_4z_\mu+a_2a_3z_\nu,z_{\mu+1},\ldots,z_{\nu-1},z_{\nu+1},\ldots,z_n\right),\\[1.5mm]
s=\left(c_1,c_2,\dots, c_{\mu-1},a_1a_4c_\mu+a_2a_3c_\nu,c_{\mu+1},\ldots,c_{\nu-1},c_{\nu+1},\ldots,c_n\right)\\[1.5mm]
 y_1=(z_1,z_2,\ldots,z_{\mu-1},z_{\mu+1},\ldots,z_n) \;\text{and}\; \;s_1=(c_1,c_2,\ldots,c_{\mu-1},c_{\mu+1},\ldots,c_n).\end{array}\eeas
\noindent 
We establish the following results, respectively, for the finite order transcendental entire functions of several complex variables satisfying the PDDEs (\ref{e1})-(\ref{e4}).
\begin{theorem}\label{TT1} Let $c=(c_1,c_2,\ldots,c_n)\in\mathbb{C}^n\setminus\{(0,0,\ldots,0)\}$ and $a_1,a_2,a_3,a_4\in\mathbb{C}\setminus\{0\}$. Let $f(z)$ be a finite order transcendental entire function on $\mathbb{C}^n$ that satisfies (\ref{e1}). Then $f(z)$ has one of the following forms:
\item[(I)] \beas f(z)=\left\{\begin{array}{l}
\frac{2K_3}{a_1\beta_\mu}e^{\frac{1}{2}\sum_{j=1}^n \beta_jz_j+\frac{1}{2}g_1(z)+\frac{1}{2}\beta}+g_2(y_1),\quad\text{when}\quad \beta_\mu\not=0;\\[2mm]
\frac{K_3 z_\mu}{a_1} e^{\frac{1}{2}\sum_{\substack{j=1,j\not=\mu}}^n \beta_jz_j+\frac{1}{2}g_1(z)+\frac{1}{2}\beta}+g_4(y_1),\quad\text{when}\quad \beta_\mu=0\end{array}\right.\eeas
with $g(z)=\sum_{j=1}^n \beta_jz_j+g_1(z)+\beta$;
\item[(II)]\beas f(z)=\left\{\begin{array}{l}
\frac{K_1e^{\sum_{j=1}^nb_jz_j+\xi_1(z)+A}}{2a_1b_\mu}+\frac{K_2e^{\sum_{j=1}^nd_jz_j+\xi_2(z)+B}}{2a_1d_\mu}+g_6(y_1),\quad\text{when}\quad b_\mu\not=0,d_\mu\not=0;\\[2mm]
\frac{K_1e^{\sum_{j=1}^nb_jz_j+\xi_1(z)+A}}{2a_1b_\mu}+\frac{K_2z_\mu e^{\sum_{j=1,j\not=\mu}^nd_jz_j+\xi_2(z)+B}}{2a_1}+g_7(y_1),\quad\text{when}\quad b_\mu\not=0, d_\mu=0;\\[2mm]
\frac{K_1z_\mu e^{\sum_{j=1,j\not=\mu}^nb_jz_j+\xi_1(z)+A}}{2a_1}+\frac{K_2e^{\sum_{j=1}^nd_jz_j+\xi_2(z)+B}}{2a_1d_\mu}+g_8(y_1),\quad\text{when}\quad b_\mu=0,d_\mu\not=0;\\[2mm]
\frac{K_1z_\mu e^{\sum_{j=1,j\not=\mu}^nb_jz_j+\xi_1(z)+A}}{2a_1}+\frac{K_2z_\mu e^{\sum_{j=1,j\not=\mu}^nd_jz_j+\xi_2(z)+B}}{2a_1}+g_9(y_1),\quad\text{when}\quad b_\mu=0=d_\mu
\end{array}\right.\eeas
with $g(z)=\sum_{j=1}^n(b_j+d_j)z_j+\xi_1(z)+\xi_2(z)+A+B$,
where $\beta_j,b_j, d_j,t_j, K_i,\beta,A,B\in\mathbb{C}$ $(1\leq j\leq n\;\text{and}\;1\leq i\leq 4)$ with $K_1K_2=1$, $K_3^2+K_4^2=1$, 
\beas &&e^{\frac{1}{2}\sum_{j=1}^n\beta_jc_j}\equiv \frac{a_1}{a_3K_3}\left(\frac{K_4\beta_\mu}{2}-\frac{a_4K_3}{4a_1}\beta_\mu^2-\frac{a_2K_3}{a_1}\right),\\
&&e^{\sum_{j=1}^nb_jc_j}\equiv -\frac{i a_1b_\mu+a_4b_\mu^2+a_2}{a_3},\; e^{\sum_{j=1}^nd_jc_j}\equiv -\frac{-i a_1d_\mu+a_4d_\mu^2+a_2}{a_3},\eeas  
$G(z)$ $(G\equiv g_1,\xi_1,\xi_2)$ is a polynomial defined in (\ref{K1}) with $G(z)\equiv 0$, when $G(z)$ contain the variable $z_\mu$ and $g_k(y_1)$ $(k=2,4,6,7,8,9)$ are finite order entire functions satisfying 
\beas \left\{\begin{array}{lll}
a_3g_2(y_1+s_1)+a_2g_2(y_1)\equiv 0\quad\text{and}\quad a_3g_6(y_1+s_1)+a_2g_6(y_1)\equiv 0,\\[2mm]
a_2g_4(y_1)+a_3g_4(y_1+s_1)\equiv\left(K_4+\frac{a_2c_\mu K_3}{a_1}\right)e^{\frac{1}{2}\sum_{\substack{j=1,j\not=\mu}}^n \beta_jz_j+\frac{1}{2}g_1(z)+\frac{1}{2}\beta},\\[2mm]
a_3 g_7(y_1+s_1)+a_2 g_7(y_1)\equiv\frac{K_2}{2}\left(\frac{a_2c_\mu}{a_1}+i\right)e^{\sum_{j=1,j\not=\mu}^nd_jz_j+\xi_2(z)+B},\\[2mm]
a_3g_8(y_1+s_1)+a_2g_8(y_1)\equiv \frac{K_1}{2}\left(\frac{a_2c_\mu}{a_1}-i\right)e^{\sum_{j=1,j\not=\mu}^nb_jz_j+\xi_1(z)+A},\\[2mm]
a_3g_9(y_1+s_1)+a_2g_9(y_1)\equiv \frac{K_1}{2}\left(\frac{a_2c_\mu}{a_1}-i\right)e^{\sum_{j=1,j\not=\mu}^nb_jz_j+\xi_1(z)+A}\\[2mm]
\qquad\qquad\qquad\qquad\qquad\quad+\frac{K_2}{2}\left(\frac{a_2c_\mu}{a_1}+i\right)e^{\sum_{j=1,j\not=\mu}^nd_jz_j+\xi_2(z)+B}.
\end{array}\right.\eeas
\end{theorem}
The following examples illustrate that the forms of the solutions presented in \textrm{Theorem \ref{TT1}} are precise.
\begin{example} Let $c=((2/5)\ln (4/15),\ln (15/4), \ln(4/15))\in\mathbb{C}^3$. It is easy to see that 
\beas f(z)
=\frac{2}{15}e^{(z_2^2+z_3^2+2z_2z_3+5z_1+7z_2+3z_3+1)/2}\eeas satisfies the PDDE
\beas \left(3\frac{\pa f(z)}{\pa z_1}\right)^2+\left(5f(z)-3f(z+c)+1\frac{\pa^2 f(z)}{\pa z_1^2}\right)^2=e^{z_2^2+z_3^2+2z_2z_3+5z_1+7z_2+3z_3+1}.\eeas
\end{example}
\begin{example} Let $c=(2\ln 3,-\ln 4,2\pi i/3)\in\mathbb{C}^3$. It is easy to see that 
\beas&& f(z)=e^{(z_1+2z_2+3z_3+5)/2}\quad\text{satisfies the PDDE}\\&& \left(2\frac{\pa f(z)}{\pa z_1}\right)^2+\left(f(z)+3f(z+c)+5\frac{\pa^2 f(z)}{\pa z_1^2}\right)^2=e^{z_1+2z_2+3z_3+5}.\eeas\end{example} 
\begin{example} Let $c=(3,-1,1)\in\mathbb{C}^3$. Then it is easy to see that 
\beas f(z)=\frac{1}{12\pi^2 i}e^{\pi i(2z_1+z_2+3z_3)+7}+\frac{1}{18\pi^2 i}e^{\pi i(3z_1+2z_2+4z_3)+5}\eeas satisfies the PDDE 
\beas \left(3\pi \frac{\pa f(z)}{\pa z_1}\right)^2+\left(5\pi^2 f(z)+5\pi^2 f(z+c)+\frac{\pa^2 f}{\pa z_1^2}\right)^2=e^{\pi i(5z_1+3z_2+7z_3)+12}.\eeas
\end{example}
\begin{theorem}\label{TT2} Let $c=(c_1,c_2,\ldots,c_n)\in\mathbb{C}^n\setminus\{(0,0,\ldots,0)\}$ and $a_i(\not=0)\in\mathbb{C}$ $(1\leq i\leq 4)$ with $a_1^2+a_3^2\not=0$. Let $f(z)$ be a finite order transcendental entire function on $\mathbb{C}^n$ that satisfies (\ref{e2}). Then $f(z)$ has one of the following forms:
\item[(I)] $f(z)=h_3(y)$, where $h_3(y)$ is a finite order transcendental entire function satisfying
\beas a_2a_3\frac{\pa h_3(y)}{\pa z_\mu}\equiv a_1a_4\frac{\pa h_3(y)}{\pa z_\nu}\;\;\text{and}\;\;a_1(h_3(y+s)-h_3(y))+a_2\frac{\pa h_3(y)}{\pa z_\mu}=\frac{a_1}{\sqrt{a_1^2+a_3^2}}e^{g(z)/2}.\eeas
\item[(II)]
\beas f(z)=\left\{\begin{array}{lll}
\frac{2\left(a_3K_3-a_1K_4\right)e^{\frac{1}{2}\sum_{j=1}^n \beta_jz_j+\frac{1}{2}g_1(z)+\frac{1}{2}\beta}}{\left(a_2a_3\beta_\mu-a_1a_4\beta_\nu\right)}+h_2(y),\\\\
\frac{\left(a_3i-a_1\right)K_1e^{\sum_{j=1}^n b_jz_j+\xi_1(z)+A}}{2i\left(a_2a_3b_\mu-a_1a_4b_\nu\right)}+\frac{\left(a_3i+a_1\right)K_2e^{\sum_{j=1}^n d_jz_j+\xi_2(z)+B}}{2i\left(a_2a_3d_\mu-a_1a_4d_\nu\right)}+h_5(y),
\end{array}\right.\eeas
where $\beta_j, b_j,d_j,\beta,A,B, K_i\in\mathbb{C}$  $(1\leq j\leq n\;\text{and}\;1\leq i\leq 4)$ with $K_1K_2=1$, $K_3^2+K_4^2=1$, $a_3K_3-a_1K_4\not=0$, $a_2a_3\chi_\mu-a_1a_4\chi_\nu\not=0$ $(\chi_i=\beta_i, b_i,d_i)$, $G(z)$ $(G\equiv g_1,\xi_1,\xi_2)$ is a polynomial defined in (\ref{K1}) with $G(z)\equiv 0$, when $G(z)$ contain the variable(s) $z_\mu$ or $z_\nu$ or both with $a_4 K_3\frac{\pa g_1(z)}{\pa z_\nu}-a_2 K_4\frac{\pa g_1(z)}{\pa z_\mu}\not=0$, $a_4 i\frac{\pa \xi_k(z)}{\pa z_\nu}+(-1)^k a_2 \frac{\pa \xi_k(z)}{\pa z_\mu}\not=0$ $(k=1,2)$ and 
\beas&& e^{\frac{1}{2}\sum_{j=1}^n\beta_jc_j}-1\equiv \frac{a_4\beta_\nu K_3-a_2\beta_\mu K_4}{a_1K_4-a_3K_3}, e^{\sum_{j=1}^nb_jc_j}-1\equiv \frac{a_4b_\nu i-a_2b_\mu}{a_1-ia_3}\\[2mm]\text{and}
&&e^{\sum_{j=1}^nd_jc_j}-1\equiv -\frac{a_4d_\nu i+a_2d_\mu}{a_1+ia_3}.\eeas 
and $h_j(y)$ $(j=2,4)$ are finite order entire functions satisfying 
\beas\left\{\begin{array}{lll}
 a_2a_3\frac{\pa h_j(y)}{\pa z_\mu}\equiv a_1a_4\frac{\pa h_j(y)}{\pa z_\nu},\\\\
h_2(y+s)-h_2(y)\equiv \frac{a_4\beta_\nu K_3-a_2\beta_\mu K_4}{a_2a_3\beta_\mu-a_1a_4\beta_\nu}e^{\frac{1}{2}\sum_{j=1}^n \beta_jz_j+\frac{1}{2}g_1(z)+\frac{1}{2}\beta}-\frac{a_2}{a_1}\frac{\pa h_2(y)}{\pa z_\mu},\\\\
h_5(y+s)-h_5(y)\equiv -\frac{a_2}{a_1}\frac{\pa h_5(y)}{\pa z_\mu}\equiv -\frac{a_4}{a_3}\frac{\pa h_5(y)}{\pa z_\nu}.\end{array}\right.\eeas
\end{theorem}
The following examples illustrate that the forms of the solutions presented in \textrm{Theorem \ref{TT2}} are precise.
\begin{example} Let $c=((1/12)\ln((27\sqrt{2}-1)/(3\sqrt{2})),\pi i/12,\pi i/2)\in\mathbb{C}^3$. Then it is easy to see that $f(z)=e^{12z_1+6z_2+z_3+9}$
satisfies the PDDE 
\beas\left(3\Delta f(z)+2\frac{\pa f(z)}{\pa z_1}\right)^2+\left(3\Delta f(z)+4\frac{\pa f(z)}{\pa z_2}\right)^2=e^{24z_1+12z_2+2z_3+18}.\eeas
\end{example}
\begin{example} Let $c=(\pi i,\pi i,2\pi i)\in\mathbb{C}^3$. Clearly 
\beas f(z)=\frac{(3i-2)e^{2z_1+3z_2+\frac{1}{2}z_3+7}}{2i(12+8i)}+\frac{(3i+2)e^{-2z_1+3z_2+\frac{3}{2}z_3+9}}{2i(-12+8i)}\eeas 
satisfies the PDDE 
\beas \left(2\Delta f(z)+2\frac{\pa f(z)}{\pa z_1}\right)^2+\left(3\Delta f(z)+\frac{4}{3i}\frac{\pa f(z)}{\pa z_2}\right)^2=e^{6z_2+2z_3+16}.\eeas\end{example}
\begin{theorem}\label{TT3} Let $c=(c_1,c_2,\ldots,c_n)\in\mathbb{C}^n\setminus\{(0,0,\ldots,0)\}$ and $a_1,a_2,a_3,a_4\in\mathbb{C}\setminus\{0\}$. 
Let $f(z)$ be a finite order transcendental entire function on $\mathbb{C}^n$ that satisfies (\ref{e3}). Then $f(z)$ has one of the following forms: 
\item[(I)] $f(z)=\pm \frac{1}{a_1}e^{\frac{g(z-c)}{2}}$, where $g(z)$ is a non-constant polynomial on $\mathbb{C}^n$ satisfying 
\beas 2a_2\frac{\pa g(z)}{\pa z_\mu}+a_3\left(\frac{\pa g(z)}{\pa z_\mu}\right)^2+2a_3\frac{\pa^2 g(z)}{\pa z_\mu^2}\equiv 0;\eeas
\item[(II)] \beas f(z)=\left\{\begin{array}{lll}
\frac{K_3}{a_1}e^{\frac{1}{2}\sum_{j=1}^n \beta_j(z_j-c_j)+\frac{1}{2}g_1(z)+\frac{1}{2}\beta}, \\\\
\frac{1}{2a_1}\left(K_1e^{\sum_{j=1}^nb_j(z_j-c_j)+\xi_1(z)+A}+K_2e^{\sum_{j=1}^nd_j(z_j-c_j)+\xi_2(z)+B}\right),
\end{array}\right.\eeas
where $\beta_j,b_j,d_j,t_j,\beta,A,B, K_i(\not=0)\in\mathbb{C}$ $(1\leq j\leq n\;\text{and}\;1\leq i\leq 4)$ with $\beta_\mu\not=0$, $b_\mu\not=0$, $d_\mu\not=0$, $K_1K_2=1$, 
$K_3^2+K_4^2=1$, $G(z)$ $(G\equiv g_1,\xi_1,\xi_2)$ is a polynomial defined in (\ref{K1}) with $G(z)\equiv 0$, when $G(z)$ contain the variable $z_\mu$, 
\beas &&e^{\frac{1}{2}\sum_{j=1}^n \beta_jc_j}\equiv \frac{K_3}{K_4}\left(\frac{a_2}{2a_1}\beta_\mu+\frac{a_3}{4a_1}\beta_\mu^2\right),e^{\sum_{j=1}^nb_jc_j}\equiv i\left(a_3b_\mu^2+a_2b_\mu\right)/a_1\\\text{and}
&&e^{\sum_{j=1}^nd_jc_j}\equiv -i\left(a_3d_\mu^2+a_2d_\mu\right)/a_1.\eeas
\end{theorem}
The following examples illustrate that the forms of the solutions presented in \textrm{Theorem \ref{TT3}} are precise.
\begin{example} Let $c=(1,2,3)\in\mathbb{C}^3$. Clearly $f(z)=e^{\frac{1}{2}(3z_1+5z_2+z_3-9)}/5$ satisfies the PDDE 
\beas 25f^2(z+c)+\left(-6\frac{\pa f(z)}{\pa z_1}+4\frac{\pa^2 f(z)}{\pa z_1^2}\right)^2=e^{3z_1+5z_2+z_3+7}.\eeas 
\end{example}
\begin{example} Let $c=(-\pi i,\pi i,3\pi i)\in\mathbb{C}^3$. It is easy to see that 
\beas f(z)=\frac{1}{4i}e^{(3z_2-z_3)^3+z_1+3z_2+2z_3+7}+\frac{1}{4i}e^{(3z_2-z_3)^2+2z_1+5z_2+z_3+5}\eeas 
satisfies the PDDE 
\beas -4f^2(z+c)+\left(5\frac{\pa f(z)}{\pa z_1}-3\frac{\pa^2 f(z)}{\pa z_1^2}\right)^2=e^{(3z_2-z_3)^3+(3z_2-z_3)^2+3z_1+8z_2+3z_3+12}.\eeas 
\end{example}
\begin{theorem}\label{TT4} Let $c=(c_1,c_2,\ldots,c_n)\in\mathbb{C}^n\setminus\{(0,0,\ldots,0)\}$ and $a_1,a_2,a_3,a_4\in\mathbb{C}\setminus\{0\}$. 
Let $f(z)$ be a finite order transcendental entire function on $\mathbb{C}^n$ that satisfies (\ref{e4}). Then $f(z)$ has one of the following forms: 
\item[(I)] $f(z)=\pm \frac{1}{a_1}e^{\frac{g(z-c)}{2}}$, where $g(z)$ is a non-constant polynomial on $\mathbb{C}^n$ satisfying 
\beas 2a_2\frac{\pa g(z)}{\pa z_\mu}+a_3\frac{\pa g(z)}{\pa z_\mu}\frac{\pa g(z)}{\pa z_\nu}+2a_3\frac{\pa^2 g(z)}{\pa z_\mu\pa z_\nu}\equiv 0;\eeas
\item[(II)] \beas f(z)=\left\{\begin{array}{lll}
\frac{K_3}{a_1}e^{\frac{1}{2}\sum_{j=1}^n \beta_j(z_j-c_j)+\frac{1}{2}g_1(z)+\frac{1}{2}\beta}, \\\\
\frac{1}{2a_1}\left(K_1e^{\sum_{j=1}^nb_j(z_j-c_j)+\xi_1(z)+A}+K_2e^{\sum_{j=1}^nd_j(z_j-c_j)+\xi_2(z)+B}\right),
\end{array}\right.\eeas
where $\beta_j,b_j,d_j,t_j,\beta,A,B, K_i(\not=0)\in\mathbb{C}$ $(1\leq j\leq n\;\text{and}\;1\leq i\leq 4)$ with $\beta_\mu\not=0$, $b_\mu\not=0$, $d_\mu\not=0$, $K_1K_2=1$, 
$K_3^2+K_4^2=1$, $G(z)$ $(G\equiv g_1,\xi_1,\xi_2)$ is a polynomial defined in (\ref{K1}) with $G(z)\equiv 0$, when $G(z)$ contain the variable $z_\mu$,
\beas &&e^{\frac{1}{2}\sum_{j=1}^n \beta_jc_j}\equiv \frac{K_3}{K_4}\left(\frac{a_2}{2a_1}\beta_\mu+\frac{a_3}{4a_1}\beta_\mu \beta_\nu\right),e^{\sum_{j=1}^nb_jc_j}\equiv i\left(a_3b_\mu b_\nu+a_2b_\mu\right)/a_1\\[2mm]\text{and}
&&e^{\sum_{j=1}^nd_jc_j}\equiv -i\left(a_3d_\mu d_\nu+a_2d_\mu\right)/a_1.\eeas
\end{theorem}
The following examples illustrate that the forms of the solutions presented in \textrm{Theorem \ref{TT4}} are precise.
\begin{example} Let $c=(\pi i,3,-\pi i)\in\mathbb{C}^3$. Clearly $f(z)=e^{\frac{1}{2}(7z_1-3z_2+5z_3+18)}/2$ satisfies the PDDE 
\beas 4f^2(z+c)+\left(6\frac{\pa f(z)}{\pa z_1}+4\frac{\pa^2 f(z)}{\pa z_1 \pa z_2}\right)^2=e^{7z_1-3z_2+5z_3+9}.\eeas 
\end{example}
\begin{example} Let $c=(-\pi i,\pi i,\pi i)\in\mathbb{C}^3$. It is easy to see that 
\beas f(z)=\frac{1}{24i}e^{2z_1+3z_2+4z_3+5}+\frac{1}{24i}e^{3z_1+z_2+5z_3+7}\eeas 
satisfies the PDDE 
\beas -144 f^2(z+c)+\left(-9\frac{\pa f(z)}{\pa z_1}+5\frac{\pa^2 f(z)}{\pa z_1\pa z_2}\right)^2=e^{5z_1+4z_2+9z_3+12}.\eeas 
\end{example}
\indent Nevanlinna's theory of several complex variables, the difference analogue of the lemma on the logarithmic derivative in several complex variables \cite{4,14}, and Lagrange's auxiliary equations \cite[Chapter 2]{430} for quasi-linear PDEs were utilized in the proof of the paper's main results.
\section {Some Lemmas} 
The following lemmas are essential to this paper and will be used to prove the main results.
\begin{lemma}\label{lem1}\cite[\textrm{Lemma 1.5}, P. 239]{12i}\cite[\textrm{Lemma 3.1, P. 211}]{12} Let $f_j\not\equiv 0\;(j=1, 2,3)$ be meromorphic functions on $\mathbb{C}^n$ such that $f_1$ is not constant and $f_1+f_2+f_3\equiv 1$ with
\beas \sum_{j=1}^3\left\{N_2(r, 0; f_j)+2\ol N(r, f_j)\right\}<\lambda T(r, f_1)+O(\log^+T(r, f_1))\eeas
holds as $r\to\infty$ out side of a possible exceptional set of finite linear measure,
where $\lambda<1$ is a positive number. Then, either $f_2\equiv 1$ or $f_3\equiv 1$. \end{lemma}
\begin{Canonical function}
Let $f(z)$ be an entire function on $\mathbb{C}^n$ ($n>1$) such that $f(0)\not=0$ and
$\rho(n(r,0,f))<\infty$. Let $q$ be the smallest integer such that the integral
\beas \int_0^\infty \frac{n(r,0,f)}{r^{q+2}}dr \;\;\text{converges}.\eeas
Then there exists an entire function $\phi(z)$ satisfying the following conditions:
\item[(i)] The function $f(z)\phi^{-1}(z)$ is an entire function on $\mathbb{C}^n$ and does not vanish.
\item[(ii)] The expansion of the function $\ln\phi(z)$ in the neighborhood of the origin has the form:
\beas \ln\phi(z)=\sum_{\Vert k\Vert=q+1}^\infty a_kz^k.\eeas
\item[(iii)] For any $R>0$, 
\bs\beas\ln M_\phi(R)\leq C_{n,q} R^q\left\{\int_0^R \frac{n(r,0,f)}{r^{q+1}}dr+R\int_R^\infty \frac{n(r,0,f)}{r^{q+2}}dr\right\}. \eeas\es
where $C_{n,q}$ is a constant and $ M_\phi(R)=\max\limits_{|z|\leq R}|\phi(z)|$. This function $\phi(z)$ is called the canonical function
(see \cite[Theorem 4.3.2, P. 245]{19}).
\end{Canonical function}
\begin{lemma}\label{lem2}\cite[Theorem 4.3.4, P. 247]{19} Let $f(z)$ be an entire function on $\mathbb{C}^n$ such that $f(0)\not=0$ and $\rho(N(r,0,f))<\infty$. 
Then there exists an entire function $g(z)$ and a canonical function $\phi(z)$
such that $f(z)=\phi(z)e^{g(z)}$.
\end{lemma}
\begin{lemma}\label{lem8}\cite[\textrm{Lemma 2.1}, P. 282]{12ai} \cite[\textrm{Lemma 3.58}]{12} If $g$ is a transcendental entire function on $\mathbb{C}^n$ and if $f$ is a meromorphic function of positive order on $\mathbb{C}$, then $f\circ g$ is of infinite order.\end{lemma}
\begin{lemma}\label{lem9}\cite[\textrm{Proposition 3.2}, P.240]{12ii} \cite[\textrm{Lemma 3.59}]{12} Let $P$ be a non-constant entire function on $\mathbb{C}^n$. Then
\beas \rho(e^P)=\left\{\begin{array}{clcr}\deg(P)&:\text{if}\; P \;\text{is a polynomial},\\
+\infty&:\text{otherwise}.\end{array}\right.\eeas
\end{lemma}
\begin{lemma}\label{lem4}\cite[\textrm{Theorem 2.1}, P. 242]{12i}\cite[\textrm{Lemma 1.106}]{12}  Suppose that $a_0(z), a_1(z),$ $\ldots,a_m(z)$ $(m\geq 1)$ are meromorphic functions on $\mathbb{C}^n$ and $g_0(z), g_1(z),
\ldots,g_m(z)$ are entire functions on $\mathbb{C}^n$ such that $g_j(z)-g_k(z)$ are not constants for $0\leq j < k \leq n$. If $\sum_{j=0}^na_j(z)e^{g_j(z)}\equiv 0$ and $T(r, a_j)=o(T(r))$, $j=0,1,\ldots,n$ hold as $r\to\infty$ out side of a possible exceptional set of finite linear measure, where $T(r)=\min\limits_{0\leq j < k \leq n} T (r, e^{g_j-g_k})$,
then $a_j(z) \equiv 0$ $(j =0,1,2,\ldots,n)$.\end{lemma}
\begin{lemma}\label{lem20}\cite[\textrm{Lemma 3.2}, P. 385]{51} Let $f$ be a non-constant meromorphic function on $\mathbb{C}^n$. Then for any $I\in \mathbb{Z}_+^n$, $T(r, \pa^I f)=O(T(r,f))$ 
for all $r$ except possibly a set of finite Lebesgue measure, where $I=\left(i_1,i_2,\ldots,i_n\right)\in\mathbb{Z}_+^n$ denotes a multiple index with $\Vert I\Vert=i_1+i_2+\cdots+i_n$, $\mathbb{Z}_+=\mathbb{N}\cup\{0\}$, and $\pa^I f=\frac{\pa^{\Vert I\Vert} f}{\pa z_1^{i_1}\cdots\pa z_n^{i_n}}$.  
\end{lemma}
\section{Proofs of the main results}  
\begin{proof}[\bf{Proof of Theorem \ref{TT1}}] Let $f$ be a finite order transcendental entire function on $\mathbb{C}^n$ that satisfies (\ref{e1}), where $g(z)$ is a non-constant polynomial on $\mathbb{C}^n$.
Now (\ref{e1}) can be written as 
\beas\prod_{j=1}^2\left(\frac{a_1\frac{\pa f(z)}{\pa z_\mu}}{e^{g(z)/2}}-(-1)^j i\frac{a_2f(z)+a_3f(z+c)+a_4\frac{\pa^2 f(z)}{\pa z_\mu^2}}{e^{g(z)/2}}\right)=1. \eeas
Here $\left(a_1\frac{\pa f(z)}{\pa z_\mu}\pm i\left(a_2f(z)+a_3f(z+c)+a_4\frac{\pa^2 f(z)}{\pa z_\mu^2}\right)\right)/e^{\frac{g(z)}{2}}$ are finite order transcendental entire functions and have no zeros on $\mathbb{C}^n$. In view of the \textrm{lemma \ref{lem2}}, we have
\beas&&\frac{a_1\frac{\pa f(z)}{\pa z_\mu}}{e^{g(z)/2}}+i\frac{a_2f(z)+a_3f(z+c)+a_4\frac{\pa^2 f(z)}{\pa z_\mu^2}}{e^{g(z)/2}}=K_1e^{P(z)}\\[2mm]\text{and}
&&\frac{a_1\frac{\pa f(z)}{\pa z_\mu}}{e^{g(z)/2}}-i\frac{a_2f(z)+a_3f(z+c)+a_4\frac{\pa^2 f(z)}{\pa z_\mu^2}}{e^{g(z)/2}}=K_2e^{-P(z)}, \eeas
where $K_1,K_2\in\mathbb{C}\setminus\{0\}$ such that $K_1K_2=1$ and $P(z)$ is an entire function on $\mathbb{C}^n$. Thus, we have
\be\label{eq1}a_1\frac{\pa f(z)}{\pa z_\mu}=\frac{K_1e^{\gamma_1(z)}+K_2e^{\gamma_2(z)}}{2}, a_2f(z)+a_3f(z+c)+a_4\frac{\pa^2 f(z)}{\pa z_\mu^2}=\frac{K_1e^{\gamma_1(z)}-K_2e^{\gamma_2(z)}}{2i},\ee
where $\gamma_1(z)=P(z)+g(z)/2$ and $\gamma_2(z)=-P(z)+g(z)/2$. Using \textrm{Lemmas \ref{lem8}}, \ref{lem9} and \ref{lem20}, it follows from (\ref{eq1}) that $P(z)$ is a 
polynomial on $\mathbb{C}^n$, as $\rho(f(z))<\infty$ and $g(z)$ is a non-constant polynomial on $\mathbb{C}^n$. 
By partial differentiation with respect to $z_\mu$ on both sides of the first equation of (\ref{eq1}), we get
\bea\label{eq3}\frac{\pa^2 f(z)}{\pa z_\mu^2}=\frac{K_1e^{\gamma_1(z)}\frac{\pa \gamma_1(z)}{\pa z_\mu}+K_2e^{\gamma_2(z)}\frac{\pa \gamma_2(z)}{\pa z_\mu}}{2a_1}.\eea  
Using (\ref{eq3}), we deduce from the second equation of (\ref{eq1}) that 
\be\label{eq4}
a_2f(z)+a_3f(z+c)=K_1e^{\gamma_1(z)}\left(\frac{1}{2i}-\frac{a_4}{2a_1}\frac{\pa \gamma_1(z)}{\pa z_\mu}\right)-K_2e^{\gamma_2(z)}\left(\frac{1}{2i}+\frac{a_4}{2a_1}\frac{\pa \gamma_2(z)}{\pa z_\mu}\right).\ee
Differentiating partially (\ref{eq4}) with respect to $z_\mu$, we get
\bea\label{eq5}&&a_2\frac{\pa f(z)}{\pa z_\mu}+a_3\frac{\pa f(z+c)}{\pa z_\mu}=K_1e^{\gamma_1(z)}\left(\frac{1}{2i}\frac{\pa \gamma_1(z)}{\pa z_\mu}-\frac{a_4}{2a_1}\left(\frac{\pa \gamma_1(z)}
{\pa z_\mu}\right)^2-\frac{a_4}{2a_1}\frac{\pa^2\gamma_1(z)}{\pa z_\mu^2}\right)\nonumber\\[2mm]
&&-K_2e^{\gamma_2(z)}\left(\frac{1}{2i}\frac{\pa \gamma_2(z)}{\pa z_\mu}+\frac{a_4}{2a_1}\left(\frac{\pa \gamma_2(z)}{\pa z_\mu}\right)^2+\frac{a_4}{2a_1}\frac{\pa^2\gamma_2(z)}{\pa z_\mu^2}\right). \eea
From (\ref{eq1}) and (\ref{eq5}), we obtain 
\bea
\label{tr}&&\Gamma_1(z)e^{\gamma_1(z)-\gamma_1(z+c)}+\Omega_1(z) e^{\gamma_2(z)-\gamma_1(z+c)}-\frac{K_2}{K_1}e^{\gamma_2(z+c)-\gamma_1(z+c)}\equiv 1, \\\text{where} 
\label{yt} &&\Gamma_1(z)=\frac{a_1}{a_3}\left(\frac{1}{i}\frac{\pa \gamma_1(z)}{\pa z_\mu}-\frac{a_4}{a_1}\left(\frac{\pa \gamma_1(z)}{\pa z_\mu}\right)^2-\frac{a_4}{a_1}\frac{\pa^2\gamma_1(z)}{\pa z_\mu^2}-\frac{a_2}{a_1}\right)\\[2mm]\text{and}
\label{yt1}&&\Omega_1(z)=-\frac{a_1K_2}{a_3K_1}\left(\frac{1}{i}\frac{\pa \gamma_2(z)}{\pa z_\mu}+\frac{a_4}{a_1}\left(\frac{\pa \gamma_2(z)}{\pa z_\mu}\right)^2+\frac{a_4}{a_1}\frac{\pa^2\gamma_2(z)}{\pa z_\mu^2}+\frac{a_2}{a_1}\right).\eea
It is necessary to consider the following cases individually.\\[2mm]
{\bf Case 1.} Let $e^{\gamma_2(z+c)-\gamma_1(z+c)}$ be constant. Then, $\gamma_2(z+c)-\gamma_1(z+c)$ must be a constant, say $k\in\mathbb{C}$. 
This implies that 
$P(z+c)\equiv-k/2$, a constant. 
From (\ref{eq1}), we have
\bea\label{eq6} a_1\frac{\pa f(z)}{\pa z_\mu}=K_3e^{\frac{g(z)}{2}}\quad\text{and}\quad a_2f(z)+a_3f(z+c)+a_4\frac{\pa^2 f(z)}{\pa z_\mu^2}=K_4e^{\frac{g(z)}{2}},\eea
where $K_3=\frac{K_1\rho+K_2\rho^{-1}}{2}$, $K_4=\frac{K_1\rho -K_2\rho^{-1}}{2i}$, $e^{-k/2}=\rho(\not=0)$ and $K_3^2+K_4^2=1$. From (\ref{eq6}), we deduce that
\bea\label{eq10}
\frac{a_1}{a_3K_3}\left(\frac{K_4}{2}\frac{\pa g(z)}{\pa z_\mu}-\frac{a_4K_3}{4a_1}\left(\frac{\pa g(z)}{\pa z_\mu}\right)^2-\frac{a_4K_3}{2a_1}\frac{\pa^2 g(z)}{\pa z_\mu^2}-\frac{a_2K_3}{a_1}\right)\equiv e^{\frac{g(z+c)-g(z)}{2}}.\eea
As $g(z)$ is a non-constant polynomial, it can be demonstrated from equation (\ref{eq10}) that $g(z+c)-g(z)$ must be a constant.
Therefore, it follows that $g(z)=\sum_{j=1}^n\beta_jz_j+g_1(z)+\beta$, 
where $\beta_i,\beta\in\mathbb{C}$ $(1\leq i\leq n)$ and $g_1(z)$ is a polynomial defined in (\ref{K1}).
From (\ref{eq10}), we have 
\be\label{n1} \frac{K_4}{2}\left(\beta_\mu+\frac{\pa g_1(z)}{\pa z_\mu}\right)-\frac{a_4K_3}{4a_1}\left(\beta_\mu+\frac{\pa g_1(z)}{\pa z_\mu}\right)^2-\frac{a_4K_3}{2a_1}\frac{\pa^2 g_1(z)}{\pa z_\mu^2}-\frac{a_2K_3}{a_1}\equiv \frac{a_3K_3}{a_1}e^{\frac{1}{2}\sum_{j=1}^n\beta_jc_j}.\ee
{\bf Sub-case 1.1.} 
If the polynomial $g_1(z)$ contains the variable $z_\mu$, a comparison of the degrees on both sides of the equation (\ref{n1}) shows that the degree of $g_1(z)$ is at most one. 
For simplicity, we still denote $g(z)=\sum_{j=1}^n\beta_jz_j+\beta$, where $\beta_j,\beta\in\mathbb{C}$ $(1\leq j\leq n)$. This implies that $g_1(z)\equiv 0$. Now we have the following cases to take into consideration.\\[2mm]
{\bf Sub-case 1.1.1.} When $\beta_\mu\not=0$. From (\ref{eq10}), we have 
\bea\label{bxw} \frac{a_1}{a_3K_3}\left(\frac{K_4\beta_\mu}{2}-\frac{a_4K_3}{4a_1}\beta_\mu^2-\frac{a_2K_3}{a_1}\right)\equiv e^{\frac{1}{2}\sum_{j=1}^n\beta_jc_j}.\eea
From (\ref{eq6}), we have
\bea\label{dse}\frac{\pa f(z)}{\pa z_\mu}=\frac{K_3}{a_1}e^{\frac{1}{2}\sum_{j=1}^n \beta_jz_j+\frac{1}{2}\beta}.\eea
The Lagrange's auxiliary equations \cite[Chapter 2]{430} of (\ref{dse}) are 
\beas \frac{dz_1}{0}=\frac{dz_2}{0}=\cdots=\frac{dz_\mu}{1}=\cdots=\frac{dz_n}{0}=\frac{df(z)}{\frac{K_3}{a_1}e^{\frac{1}{2}\sum_{j=1}^n \beta_jz_j+\frac{1}{2}\beta}}.\eeas
Note that $\alpha_j=z_j$ for $1\leq j(\not=\mu)\leq n$ and 
\beas &&df(z)=\frac{K_3}{a_1}e^{\frac{1}{2}\sum_{j=1}^n \beta_jz_j+\frac{1}{2}\beta}dz_\mu=\frac{K_3}{a_1}e^{\frac{1}{2}\beta_\mu z_\mu+\frac{1}{2}\sum_{\substack{j=1,j\not=\mu}}^n \beta_j\alpha_j+\frac{1}{2}\beta}dz_\mu,\\[2mm]\text{{\it i.e.,}}\quad &&f(z)=\frac{2K_3}{a_1\beta_\mu}e^{\frac{1}{2}\sum_{j=1}^n \beta_jz_j+\frac{1}{2}\beta}+\alpha_\mu.\eeas 
Note that after integration with respect to 
$z_\mu$, replacing $\alpha_j$ by $z_j$ for $1\leq j(\not=\mu)\leq n$, where $\alpha_j\in\mathbb{C}$ for $1\leq j\leq n$. Hence the solution is $\Phi(\alpha_1,\alpha_2,\ldots,\alpha_\mu,\ldots, \alpha_n )=0$. For simplicity, we suppose
\bea\label{n2} f(z)=\frac{2K_3}{a_1\beta_\mu}e^{\frac{1}{2}\sum_{j=1}^n \beta_jz_j+\frac{1}{2}\beta}+g_2(y_1),\eea
where $g_2(y_1)$ is a finite order entire function of $z_1,z_2,\ldots,z_{\mu-1},z_{\mu+1},\ldots,z_n$. 
From the second equation of (\ref{eq6}) with the help of (\ref{bxw}) and (\ref{n2}), we get $a_3g_2(y_1+s_1)+a_2g_2(y_1)\equiv 0$. \\[2mm]
{\bf Sub-case 1.1.2.} When $\beta_\mu=0$. 
From (\ref{eq6}) and (\ref{eq10}), we have 
\bea\label{pq1}\frac{\pa f(z)}{\pa z_\mu}=\frac{K_3}{a_1}e^{\frac{1}{2}\sum_{\substack{j=1,j\not=\mu}}^n\beta_jz_j+\frac{1}{2}\beta}\quad\textrm{and}\quad e^{\frac{1}{2}\sum_{j=1,j\not=\mu}^n\beta_jc_j}\equiv -\frac{a_2}{a_3}.\eea
Using arguments similar to those in \textrm{Sub-case 1.1.1}, we derive from (\ref{pq1}) that
\bea\label{pq2} f(z)=\frac{K_3 z_\mu}{a_1}e^{\frac{1}{2}\sum_{j=1,j\not=\mu}^n \beta_jz_j+\frac{1}{2}\beta}+g_4(y_1),\eea
where $g_4(y_1)$ is a finite order entire function of $z_1,z_2,\ldots,z_{\mu-1},z_{\mu+1},\ldots,z_n$. 
From the second equation of (\ref{eq6}) with the help of (\ref{pq1}) and (\ref{pq2}), we get
\beas a_3g_4(y_1+s_1)+a_2g_4(y_1)\equiv\left(K_4+\frac{a_2K_3c_\mu}{a_1}\right)e^{\frac{1}{2}\sum_{\substack{j=1,j\not=\mu}}^n \beta_jz_j+\frac{1}{2}\beta}.\eeas
{\bf Sub-case 1.2.} 
If $g_1(z)$ is independent of $z_\mu$, then from (\ref{n1}), we get (\ref{bxw}) and 
\beas g(z)=\sum_{j=1}^n\beta_jz_j+g_1(z)+\beta,\eeas 
where $\beta_i,\beta\in\mathbb{C}$ ($1\leq i\leq n$) and $g_1(z)$ is a polynomial defined in (\ref{K1}).\\[2mm]
{\bf Sub-case 1.2.1.} 
If $\beta_\mu\not=0$, using similar arguments as in \textrm{Sub-case 1.1.1}, we get from the first equation (\ref{eq6})
\bea\label{n3} f(z)=\frac{2K_3}{a_1\beta_\mu}e^{\frac{1}{2}\sum_{j=1}^n \beta_jz_j+\frac{1}{2}g_1\left(z\right)+\frac{1}{2}\beta}+g_3(y_1),\eea
where $g_3(y_1)$ is a finite order entire function. Using (\ref{bxw}) and (\ref{n3}), we get from the second equation of (\ref{eq6}) that 
$a_3g_3(y_1+s_1)+a_2g_3(y_1)\equiv 0$. \\[2mm]
{\bf Sub-case 1.2.2.} If $\beta_\mu=0$, then from the first equation of (\ref{eq6}), using similar arguments as in \textrm{Sub-case 1.1.2}, we have
\bea\label{n4}\frac{\pa f(z)}{\pa z_\mu}=\frac{K_3}{a_1}e^{\frac{1}{2}\sum_{\substack{j=1,j\not=\mu}}^n \beta_jz_j+\frac{1}{2}g_1\left(z\right)+\frac{1}{2}\beta}\quad\text{and}\quad e^{\frac{1}{2}\sum_{\substack{j=1,j\not=\mu}}^n\beta_jc_j}\equiv -\frac{a_2}{a_3}.\eea
We use arguments similar to those in \textrm{Sub-case 1.1.2} and deduce from (\ref{n4}) that
\beas\label{n5} f(z)=\frac{K_3}{a_1}e^{\frac{1}{2}\sum_{\substack{j=1,j\not=\mu}}^n \beta_jz_j+\frac{1}{2}g_1\left(z\right)+\frac{1}{2}\beta}z_\mu+g_5(y_1),\eeas
where $g_5(y_1)$ is a finite order entire function satisfying
\beas a_2g_5(y_1)+a_3g_5(y_1+s_1)\equiv\left(K_4+\frac{a_2K_3c_\mu}{a_1}\right)e^{\frac{1}{2}\sum_{\substack{j=1,j\not=\mu}}^n \beta_jz_j+\frac{1}{2}g_1\left(z\right)+\frac{1}{2}\beta}.\eeas
{\bf Case 2.} Let $e^{\gamma_2(z+c)-\gamma_1(z+c)}$ be non-constant. It is evident from (\ref{tr}) that $\Gamma_1(z)$ and $\Omega_1(z)$ are not simultaneously identically zero. 
Otherwise, we arrive at a contradiction. 
Let $\Gamma_1(z)\not\equiv 0$ and $\Omega_1(z)\equiv 0$. Then from (\ref{tr}), we have 
\be\label{tr1}\Gamma_1(z)e^{\gamma_1(z)-\gamma_1(z+c)}-\frac{K_2}{K_1}e^{\gamma_2(z+c)-\gamma_1(z+c)}\equiv 1,\;\text{{\it i.e.,}}\; \Gamma_1(z)e^{\gamma_1(z)}-\frac{K_2}{K_1}e^{\gamma_2(z+c)}-e^{\gamma_1(z+c)}\equiv 0.\ee
From (\ref{tr1}), it is clear that $\gamma_1(z)-\gamma_1(z+c)$ is not a constant. We claim that $\gamma_2(z+c)-\gamma_1(z)$ is non-constant. If not, let $\gamma_2(z+c)-\gamma_1(z)\equiv k$ which implies that $\gamma_2(z+c)\equiv \gamma_1(z)+k$, where $k\in\mathbb{C}$. From (\ref{tr1}), we have $\left(M_1(z)-K_2e^{k}/K_1\right)e^{\gamma_1(z)-\gamma_1(z+c)}\equiv 1$, which contradicts the fact that $\gamma_1(z)-\gamma_1(z+c)$ is not a constant. In view of \textrm{Lemma \ref{lem4}}, we get a contradiction from (\ref{tr1}). Similarly, we get a contradiction when 
$M_1(z)\equiv 0$, $N_1(z)\not\equiv 0$. Hence $M_1(z)\not\equiv 0$ and $N_1(z)\not\equiv 0$. Since $e^{\gamma_2(z+c)-\gamma_1(z+c)}$ is non-constant and it is evident that 
\beas&& N\left(r, \Gamma_1(z)e^{\gamma_1(z)-\gamma_1(z+c)}\right)= N\left(r, \Omega_1(z)e^{\gamma_2(z)-\gamma_1(z+c)}\right)=N\left(r,-K_2e^{\gamma_2(z+c)-\gamma_1(z+c)}/K_1\right)\\[2mm]
&&=N\left(r, 0; \Gamma_1(z)e^{\gamma_1(z)-\gamma_1(z+c)}\right)=N\left(r, 0; \Omega_1(z)e^{\gamma_2(z)-\gamma_1(z+c)}\right)\\[2mm]
&&=N\left(r,0; -K_2e^{\gamma_2(z+c)-\gamma_1(z+c)}/K_1\right)=S\left(r, -K_2e^{\gamma_2(z+c)-\gamma_1(z+c)}/K_1\right).\eeas
In the light of \textrm{Lemma \ref{lem1}}, it follows from (\ref{tr}) that either 
\beas \Gamma_1(z)e^{\gamma_1(z)-\gamma_1(z+c)}\equiv 1\quad\text{or}\quad \Omega_1(z) e^{\gamma_2(z)-\gamma_1(z+c)}\equiv 1.\eeas
{\bf Sub-case 2.1.} Let $\Gamma_1(z)e^{\gamma_1(z)-\gamma_1(z+c)}\equiv 1$. From (\ref{tr}), we have 
\bea\label{n6} \frac{K_1}{K_2}\Omega_1(z) e^{\gamma_2(z)-\gamma_2(z+c)}\equiv 1,\eea
where $\Gamma_1(z)$, $\Omega_1(z)$ are given in (\ref{yt}) and (\ref{yt1}) respectively. Therefore $\gamma_1(z)-\gamma_1(z+c)$ and $\gamma_2(z)-\gamma_2(z+c)$ are both constants. By means of arguments similar to those presented in \textrm{Case 1}, we have
$\gamma_1(z)=\sum_{j=1}^nb_jz_j+\xi_1(z)+A$ and $\gamma_2(z)=\sum_{j=1}^nd_jz_j+\xi_2(z)+B$, where $b_j,d_j,A,B\in\mathbb{C}$ $(1\leq j\leq n)$ and  
$\xi_k(z)$ $(k=1,2)$ are polynomials defined in (\ref{K1}). 
Therefore, we have 
\bea\label{n7}\begin{array}{lll}
 \frac{a_1}{a_3}\left(\frac{1}{i}\left(b_\mu+\frac{\pa \xi_1(z)}{\pa z_\mu}\right)-\frac{a_4}{a_1}\left(b_\mu+\frac{\pa \xi_1(z)}{\pa z_\mu}\right)^2-\frac{a_4}{a_1}\frac{\pa^2 \xi_1(z)}{\pa z_\mu^2}-\frac{a_2}{a_1}\right)\equiv e^{\sum_{j=1}^nb_jc_j}\\[2mm]
 -\frac{a_1}{a_3}\left(\frac{1}{i}\left(d_\mu+\frac{\pa \xi_2(z)}{\pa z_\mu}\right)+\frac{a_4}{a_1}\left(d_\mu+\frac{\pa \xi_2(z)}{\pa z_\mu}\right)^2+\frac{a_4}{a_1}\frac{\pa^2 \xi_1(z)}{\pa z_\mu^2}+\frac{a_2}{a_1}\right)\equiv e^{\sum_{j=1}^nd_jc_j}.\end{array}\eea
The following cases are to be considered separately.\\[2mm]
{\bf Sub-case 2.1.1.} When $\xi_k(z)$ is dependent of the variable $z_\mu$, then by comparing the degrees on both sides of (\ref{n7}), we get that $\deg(\xi_k(z))\leq 1$ for $k=1,2$. For simplicity, we still denote $\gamma_1(z)=\sum_{j=1}^nb_jz_j+A$ and $\gamma_2(z)=\sum_{j=1}^nd_jz_j+B$, where $b_j, d_j, A, B\in\mathbb{C}$ $(1\leq j\leq n)$. This implies that 
$\xi_k(z)\equiv 0$. Since $\gamma_2(z+c)-\gamma_1(z+c)$ is a non-constant polynomial, so we must have $b_j\not=d_j$ for some $j$. Now we have the following cases to take into consideration.\\[2mm]
{\bf Sub-case 2.1.1.1.} When $b_\mu\not=0$, $d_\mu\not=0$. Then, we have
\be\label{yt2} \frac{a_1}{a_3}\left(\frac{1}{i}b_\mu-\frac{a_4}{a_1}b_\mu^2-\frac{a_2}{a_1}\right)\equiv e^{\sum_{j=1}^nb_jc_j}\quad \text{and}\quad -\frac{a_1}{a_3}\left(\frac{1}{i}d_\mu+\frac{a_4}{a_1}d_\mu^2+\frac{a_2}{a_1}\right)\equiv e^{\sum_{j=1}^nd_jc_j}.\ee
We deduce from (\ref{eq1}), using similar arguments as in \textrm{Sub-case 1.1.1}, that
\be\label{5tr5}
f(z)=\frac{K_1e^{\sum_{j=1}^nb_jz_j+A}}{2a_1b_\mu}+\frac{K_2e^{\sum_{j=1}^nd_jz_j+B}}{2a_1d_\mu}+g_6(y_1), \ee
where $b_j, d_j,A, B\in\mathbb{C}$ $(1\leq j\leq n)$ and $g_6(y_1)$ is a finite order entire function of $z_1,z_2,\ldots,z_{\mu-1},\\z_{\mu+1},\ldots,z_n$. Using (\ref{yt2}) and (\ref{5tr5}), we deduce from the second equation of (\ref{eq1}) that $
a_3 g_6(y_1+s_1)+a_2g_6(y_1)\equiv 0$. \\[2mm]
{\bf Sub-case 2.1.1.2.} When $b_\mu\not=0$, $d_\mu=0$. Then, we have 
\bea\label{n8} e^{\sum_{j=1}^nb_jc_j}\equiv\frac{a_1}{a_3}\left(\frac{1}{i}b_\mu-\frac{a_4}{a_1}b_\mu^2-\frac{a_2}{a_1}\right)\quad\text{and}\quad e^{\sum_{j=1,j\not=\mu}^nd_jc_j}\equiv -\frac{a_2}{a_3}.\eea 
From (\ref{eq1}) we deduce, by means of arguments similar to those in \textrm{Sub-case 1.1.1}, that
\beas f(z)=\frac{K_1e^{\sum_{j=1}^nb_jz_j+A}}{2a_1b_\mu}+\frac{K_2 z_\mu}{2a_1}e^{\sum_{j=1,j\not=\mu}^nd_jz_j+B}+g_7(y_1),\eeas
where $b_j,d_k,A,B\in\mathbb{C}$ ($1\leq j\leq n$ and $1\leq k(\not=\mu)\leq n$) and $g_7(y_1)$ is a finite order entire function satisfying
\beas a_3 g_7(y_1+s_1)+a_2g_7(y_1)\equiv\frac{K_2}{2}\left(\frac{a_2c_\mu}{a_1}+i\right)e^{\sum_{j=1,j\not=\mu}^nd_jz_j+B}.\eeas
{\bf Sub-case 2.1.1.3.} When $b_\mu=0$, $d_\mu\not=0$. Using arguments similar to those presented in \textrm{Sub-case 2.1.1.2}, we deduce that
\beas f(z)=\frac{K_1z_\mu}{2a_1}e^{\sum_{j=1,j\not=\mu}^nb_jz_j+A}+\frac{K_2e^{\sum_{j=1}^nd_jz_j+B}}{2a_1d_\mu}+g_8(y_1),\eeas 
where $b_j,d_k, A,B\in\mathbb{C}$ ($1\leq j(\not=\mu)\leq n$ and $1\leq k\leq n$) and $g_8(y_1)$ is a finite order entire function satisfying 
\bea&& a_3g_8(y_1+s_1)+a_2g_8(y_1)\equiv \frac{K_1}{2}\left(\frac{a_2c_\mu}{a_1}-i\right)e^{\sum_{j=1,j\not=\mu}^nb_jz_j+A},\nonumber\\[2mm]
\label{n9}&&e^{\sum_{j=1,j\not=\mu}^nb_jc_j}\equiv -\frac{a_2}{a_3}\quad\text{and}\quad e^{\sum_{j=1}^nd_jc_j}\equiv -\frac{a_1}{a_3}\left(\frac{1}{i}d_\mu+\frac{a_4}{a_1}d_\mu^2+\frac{a_2}{a_1}\right).\eea 
{\bf Sub-case 2.1.1.4.} When $b_\mu=0$, $d_\mu=0$. Using arguments similar to those presented in \textrm{Sub-case 2.1.1.2}, we deduce that
\beas f(z)=\frac{K_1 z_\mu}{2a_1}e^{\sum_{j=1,j\not=\mu}^nb_jz_j+A}+\frac{K_2z_\mu}{2a_1}e^{\sum_{j=1,j\not=\mu}^nd_jz_j+B}+g_9(y_1)\eeas
where $b_j,d_j, A,B\in\mathbb{C}$ ($1\leq j(\not=\mu)\leq n$), $e^{\sum_{j=1,j\not=\mu}^nb_jc_j}\equiv -a_2/a_3\equiv e^{\sum_{j=1,j\not=\mu}^nd_jc_j}$ and $g_9(y_1)$ is a finite order entire function satisfying
\beas a_3g_9(y_1+s_1)+a_2g_9(y_1)\equiv \frac{K_1}{2}\left(\frac{a_2c_\mu}{a_1}-i\right)e^{\sum_{j=1,j\not=\mu}^nb_jz_j+A}+\frac{K_2}{2}\left(\frac{a_2c_\mu}{a_1}+i\right)e^{\sum_{j=1,j\not=\mu}^nd_jz_j+B}.\eeas
{\bf Sub-case 2.1.2.} When $\xi_k(z)$ $(k=1,2)$ is independent of the variable $z_\mu$, then from (\ref{n7}), we get (\ref{yt2}). Thus, we have 
$$ \gamma_1(z)=\sum_{j=1}^nb_jz_j+\xi_1\left(z\right)+A\quad\text{and}\quad \gamma_2(z)=\sum_{j=1}^nd_jz_j+\xi_2\left(z\right)+B,$$
where $b_j,d_j,A,B\in\mathbb{C}$ $(1\leq j\leq n)$ and $\xi_k$ $(k=1,2)$ is a polynomial defined in (\ref{K1}).\\[2mm]
{\bf Sub-case 2.1.2.1.} When $b_\mu\not=0$, $d_\mu\not=0$. Then, by using similar arguments as in \textrm{Sub-case 2.1.1.1}, that
\beas f(z)=\frac{K_1e^{\sum_{j=1}^nb_jz_j+\xi_1\left(z\right)+A}}{2a_1b_\mu}+\frac{K_2e^{\sum_{j=1}^nd_jz_j+\xi_2\left(z\right)+B}}{2a_1d_\mu}+h_1(y_1), \eeas
where $b_j,d_j,A,B\in\mathbb{C}$ $(1\leq j\leq n)$ with (\ref{yt2}) and $h_1(y_1)$ is a finite order entire function satisfying $a_3 h_1(y_1+s_1)+a_2 h_1(y_1)\equiv 0$. \\[2mm]
{\bf Sub-case 2.1.2.2.} When $b_\mu\not=0$, $d_\mu=0$. By means of arguments similar to those in \textrm{Sub-case 2.1.1.2}, that
\beas f(z)=\frac{K_1e^{\sum_{j=1}^nb_jz_j+\xi_1\left(z\right)+A}}{2a_1b_\mu}+\frac{K_2z_\mu}{2a_1}e^{\sum_{j=1,j\not=\mu}^nd_jz_j+\xi_2\left(z\right)+B}+h_2(y_1),\eeas
where $b_j,d_k,A,B\in\mathbb{C}$ ($1\leq j\leq n$ and $1\leq k(\not=\mu)\leq n$) with (\ref{n8}) and $h_2(y_1)$ is a finite order entire function satisfying 
\beas a_3 h_2(y_1+s_1)+a_2 h_2(y_1)\equiv\frac{K_2}{2}\left(\frac{a_2c_\mu}{a_1}+i\right)e^{\sum_{j=1,j\not=\mu}^nd_jz_j+\xi_2\left(z\right)+B},\\
e^{\sum_{j=1}^nb_jc_j}\equiv\frac{a_1}{a_3}\left(\frac{1}{i}b_\mu-\frac{a_4}{a_1}b_\mu^2-\frac{a_2}{a_1}\right)\quad\text{and}\quad e^{\sum_{j=1,j\not=\mu}^nd_jc_j}\equiv -\frac{a_2}{a_3}.\eeas
{\bf Sub-case 2.1.2.3.} When $b_\mu=0$, $d_\mu\not=0$. Using arguments similar to those presented in \textrm{Sub-case 2.1.1.3}, we deduce that
\beas f(z)=\frac{K_1z_\mu}{2a_1}e^{\sum_{j=1,j\not=\mu}^nb_jz_j+\xi_1\left(z\right)+A}+\frac{K_2e^{\sum_{j=1}^nd_jz_j+\xi_2\left(z\right)+B}}{2a_1d_\mu}+h_3(y_1),\eeas 
where $b_j,d_k,A,B\in\mathbb{C}$ ($1\leq j(\not=\mu)\leq n$ and $1\leq k\leq n$) with (\ref{n9}) and $h_3(y_1)$ is a finite order entire function satisfying 
\beas&& a_3h_3(y_1+s_1)+a_2h_3(y_1)\equiv \frac{K_1}{2}\left(\frac{a_2c_\mu}{a_1}-i\right)e^{\sum_{j=1,j\not=\mu}^nb_jz_j+\xi_1\left(z\right)+A}.\eeas 
{\bf Sub-case 2.1.2.4.} When $b_\mu=0$, $d_\mu=0$. Using arguments similar to those presented in \textrm{Sub-case 2.1.1.4}, we deduce that
\beas f(z)=\frac{K_1 z_\mu}{2a_1}e^{\sum_{j=1,j\not=\mu}^nb_jz_j+\xi_1\left(z\right)+A}+\frac{K_2 z_\mu}{2a_1}e^{\sum_{j=1,j\not=\mu}^nd_jz_j+\xi_2\left(z\right)+B}+h_4(y_1)\eeas
where $b_j,d_j,A,B\in\mathbb{C}$ ($1\leq j(\not=\mu)\leq n$) and $h_4(y_1)$ is a finite order entire function satisfying
\beas a_3h_4(y_1+s_1)+a_2h_4(y_1)&\equiv& \frac{K_1}{2}\left(\frac{a_2c_\mu}{a_1}-i\right)e^{\sum_{j=1,j\not=\mu}^nb_jz_j+\xi_1\left(z\right)+A}\\
&&+\frac{K_2}{2}\left(\frac{a_2c_\mu}{a_1}+i\right)e^{\sum_{j=1,j\not=\mu}^nd_jz_j+\xi_2\left(z\right)+B}\\\text{and}\qquad
e^{\sum_{j=1,j\not=\mu}^nb_jc_j}\equiv -a_2/a_3&\equiv& e^{\sum_{j=1,j\not=\mu}^nd_jc_j}.\eeas
{\bf Sub-case 2.2.} Let $\Omega_1(z) e^{\gamma_2(z)-\gamma_1(z+c)}\equiv 1$. From (\ref{tr}), we have 
\bea\label{yt3} \frac{K_1}{K_2}\Gamma_1(z) e^{\gamma_1(z)-\gamma_2(z+c)}\equiv 1,\eea
where $\Gamma_1(z)$, $\Omega_1(z)$ are given in (\ref{yt}) and (\ref{yt1}) respectively. Therefore $\gamma_2(z)-\gamma_1(z+c)$ and $\gamma_1(z)-\gamma_2(z+c)$ are both constants, say $\chi_1$ and $\chi_2$ 
respectively, where $\chi_1,\chi_2\in\mathbb{C}$. Now $\gamma_1(z)-\gamma_1(z+2c)=\left(\gamma_1(z)-\gamma_2(z+c)\right)+\left(\gamma_2(z+c)-\gamma_1(z+2c)\right)\equiv\chi_1+\chi_2$ and 
$\gamma_2(z)-\gamma_2(z+2c)=\left(\gamma_2(z)-\gamma_1(z+c)\right)+\left(\gamma_1(z+c)-\gamma_2(z+2c)\right)\equiv\chi_1+\chi_2$.
Using arguments similar to \textrm{Case 1}, we derive that
$\gamma_1(z)=\sum_{j=1}^nb_jz_j+\xi_3(z)+A$ and $\gamma_2(z)=\sum_{j=1}^nd_jz_j+\xi_4(z)+B$, where $b_j,d_j,A,B\in\mathbb{C}$ ($1\leq j\leq n$) and 
$\xi_k(z)$ $(k=3,4)$ is a polynomial defined in (\ref{K1}). Since $\gamma_2(z)-\gamma_1(z+c)$ and $\gamma_1(z)-\gamma_2(z+c)$ are both constants, thus we conclude that $b_j=d_j$ for $1\leq j\leq n$ and $\xi_3(z)\equiv \xi_4(z)$.
Therefore $\gamma_1(z+c)-\gamma_2(z+c)=A-B$, which is a constant that contradicts the fact that $e^{\gamma_2(z+c)-\gamma_1(z+c)}$ is not a constant.
This completes the proof.
\end{proof}
\begin{proof}[\bf{Proof of Theorem \ref{TT2}}] Let $f$ be a finite order transcendental entire function on $\mathbb{C}^n$ that satisfies (\ref{e2}), where $g(z)$ is a non-constant polynomial on $\mathbb{C}^n$. Now (\ref{e2}) can be written as 
\beas\prod_{j=1}^2\left(\frac{a_1\Delta f(z)+a_2\frac{\pa f}{\pa z_\mu}}{e^{\frac{g(z)}{2}}}-(-1)^j i\frac{a_3\Delta f(z)+a_4\frac{\pa f}{\pa z_\nu}}{e^{\frac{g(z)}{2}}}\right)=1. \eeas
Using arguments similar to those in \textrm{Theorem \ref{TT1}}, we obtain 
\be\label{po1}a_1\Delta f(z)+a_2\frac{\pa f(z)}{\pa z_\mu}=\frac{K_1e^{\gamma_1(z)}+K_2e^{\gamma_2(z)}}{2},\;
a_3\Delta f(z)+a_4\frac{\pa f(z)}{\pa z_\nu}=\frac{K_1e^{\gamma_1(z)}-K_2e^{\gamma_2(z)}}{2i},\ee
where $K_1,K_2\in\mathbb{C}\setminus\{0\}$ such that $K_1K_2=1$, $\gamma_1(z)=P(z)+g(z)/2$, $\gamma_2(z)=-P(z)+g(z)/2$ and
$P(z)$ is a polynomial on $\mathbb{C}^n$.
From (\ref{po1}), we deduce that 
\bea\label{po3} a_2a_3\frac{\pa f(z)}{\pa z_\mu}-a_1a_4\frac{\pa f(z)}{\pa z_\nu}
=\left(\frac{a_3}{2}-\frac{a_1}{2i}\right)K_1e^{\gamma_1(z)}+\left(\frac{a_3}{2}+\frac{a_1}{2i}\right)K_2e^{\gamma_2(z)}.\eea
Note that $\frac{\pa^2 }{\pa z_\mu\pa z_\nu}f(z)=\frac{\pa^2}{\pa z_\mu\pa z_\mu}f(z)$. 
By partially differentiating the first and second equations of (\ref{po1}) with respect to $z_\nu$ and $z_\mu$ respectively, we have
\bea\label{po4}&&a_1\frac{\pa\Delta f(z)}{\pa z_\nu}+a_2\frac{\pa^2 f(z)}{\pa z_\nu\pa z_\mu}=\frac{1}{2}\left(K_1e^{\gamma_1(z)}\frac{\pa\gamma_1(z)}{\pa z_\nu}+K_2e^{\gamma_2(z)}\frac{\pa\gamma_2(z)}{\pa z_\nu}\right)\\[2mm]\text{and}
\label{po5}&&a_3\frac{\pa\Delta f(z)}{\pa z_\mu}+a_4\frac{\pa^2 f(z)}{\pa z_\mu\pa z_\nu}=\frac{1}{2i}\left(K_1e^{\gamma_1(z)}\frac{\pa\gamma_1(z)}{\pa z_\mu}-K_2e^{\gamma_2(z)}\frac{\pa\gamma_2(z)}{\pa z_\nu}\right).\eea
From (\ref{po4}) and (\ref{po5}), we deduce that
\bea\label{po6}&&a_1a_4\frac{\pa\Delta f(z)}{\pa z_\nu}-a_2a_3\frac{\pa\Delta f(z)}{\pa z_\mu}
=\left(\frac{a_4}{2}\frac{\pa\gamma_1(z)}{\pa z_\nu}-\frac{a_2}{2i}\frac{\pa\gamma_1(z)}{\pa z_\mu}\right)K_1e^{\gamma_1(z)}\nonumber\\
&&+\left(\frac{a_4}{2}\frac{\pa\gamma_2(z)}{\pa z_\nu}+\frac{a_2}{2i}\frac{\pa\gamma_2(z)}{\pa z_\mu}\right)K_2e^{\gamma_2(z)}.\eea
Using (\ref{po3}), we get from (\ref{po6}) that
\bea\label{nb}
\Gamma_2(z) e^{\gamma_1(z)-\gamma_1(z+c)}+\frac{\left(a_1+ia_3\right)K_2}{\left(a_1-i a_3\right)K_1}e^{\gamma_2(z+c)-\gamma_1(z+c)}+\Omega_2(z)e^{\gamma_2(z)-\gamma_1(z+c)}\equiv 1,\eea
where 
\beas&&\Gamma_2(z)=(i a_4\frac{\pa\gamma_1(z)}{\pa z_\nu}-a_2\frac{\pa\gamma_1(z)}{\pa z_\mu}+\left(a_1-ia_3\right))/(a_1-ia_3)\\\text{and}
&&\Omega_2(z)=\left(i a_4\frac{\pa\gamma_2(z)}{\pa z_\nu}+a_2\frac{\pa\gamma_2(z)}{\pa z_\mu}-\left(a_1+i a_3\right)\right)K_2/(\left(a_1-i a_3\right)K_1).\eeas
The following cases occur separately.\\[2mm]
{\bf Case 1.} If $e^{\gamma_2(z+c)-\gamma_1(z+c)}$ is constant. By means of arguments similar to those in \textrm{Case 1} of \textrm{Theorem \ref{TT1}},
we deduce that $P(z+c)\equiv-k/2$, where $k\in\mathbb{C}$.
From (\ref{po1}), we have 
\bea\label{mn}a_1\Delta f(z)+a_2\frac{\pa f(z)}{\pa z_\mu}=K_3e^{\frac{g(z)}{2}}\quad
\text{and}\quad a_3\Delta f(z)+a_4\frac{\pa f(z)}{\pa z_\nu}=K_4e^{\frac{g(z)}{2}},\eea
where $K_3=(K_1\rho +K_2\rho^{-1})/2$, $K_4=(K_1\rho -K_2\rho^{-1})/(2i)$, $e^{-\frac{k}{2}}=\rho(\not=0)$ and $K_3^2+K_4^2=1$. 
From (\ref{mn}), we have
\bea\label{mn0}a_2a_3\frac{\pa f(z)}{\pa z_\mu}-a_1a_4\frac{\pa f(z)}{\pa z_\nu}=\left(a_3K_3-a_1K_4\right)e^{\frac{g(z)}{2}}.\eea
{\bf Sub-case 1.1.} When $a_3K_3-a_1K_4\not=0$.
From (\ref{mn}) and (\ref{mn0}), we deduce that
\bea\label{mn5}(a_1K_4-a_3K_3) \left(e^{\frac{g(z+c)-g(z)}{2}}-1\right)\equiv a_4K_3\frac{\pa g(z)}{\pa z_\nu}-a_2K_4\frac{\pa g(z)}{\pa z_\mu}.\eea
Since $g(z)$ is a non-constant polynomial, it is follows from (\ref{mn5}) that $g(z+c)-g(z)$ must be a constant. 
Therefore, we conclude that
$g(z)=\sum_{j=1}^n \beta_jz_j+g_1(z)+\beta$, where $\beta_j,\beta\in\mathbb{C}$ $(1\leq j\leq n)$ and $g_1(z)$ is a polynomial defined in (\ref{K1}). From (\ref{mn5}), we have 
\be\label{poo} (a_1K_4-a_3K_3)\left(e^{\frac{1}{2}\sum_{j=1}^n\beta_jc_j}-1\right)
\equiv a_4\beta_\nu K_3-a_2\beta_\mu K_4+\left(a_4K_3\frac{\pa g_1(z)}{\pa z_\nu}-a_2K_4\frac{\pa g_1(z)}{\pa z_\mu}\right).\ee
{\bf Sub-case 1.1.1.}
Let $g_1(z)$ contains the variables $z_\mu$ or $z_\nu$ or both with $a_4K_3\frac{\pa }{\pa z_\nu}g_1(z)\not= a_2K_4\frac{\pa }{\pa z_\mu}g_1(z)$, then by comparing the degrees on both sides of (\ref{poo}), we get that $\deg(g_1(z))\leq 1$. For simplicity, we still denote 
$g(z)=\sum_{j=1}^n\beta_jz_j+\beta$, where $\beta_j,\beta\in\mathbb{C}$ $(1\leq j\leq n)$. This implies that $g_1(z)\equiv 0$. Thus, we have  
\bea\label{K2}e^{\frac{1}{2}\sum_{j=1}^n\beta_jc_j}-1\equiv(a_4\beta_\nu K_3-a_2\beta_\mu K_4)/(a_1K_4-a_3K_3).\eea
The Lagrange's auxiliary equations \cite[Chapter 2]{430} of (\ref{mn0}) are 
\beas \frac{dz_1}{0}=\frac{dz_1}{0}=\cdots=\frac{dz_\mu}{a_2a_3}=\cdots=\frac{dz_\nu}{-a_1a_4}=\cdots=\frac{dz_n}{0}=\frac{df(z)}{\left(a_3K_3-a_1K_4\right)e^{\frac{1}{2}\sum_{j=1}^n \beta_jz_j+\frac{1}{2}\beta}}.\eeas
Note that $\alpha_\nu=a_1a_4z_\mu+a_2a_3z_\nu$, $\alpha_j=z_j$ ($1\leq j(\not=\mu,\nu)\leq n$) and 
\beas df(z)&=&\frac{\left(a_3K_3-a_1K_4\right)e^{\frac{1}{2}\sum_{j=1}^n \beta_jz_j+\frac{1}{2}\beta}}{a_2a_3}\;dz_1\\
&=&\frac{\left(a_3K_3-a_1K_4\right)e^{\frac{1}{2}\left(\sum_{j=1,j\not=\mu,\nu}^n\beta_j \alpha_j+\beta_\mu z_\mu+\beta_\nu\left(\frac{\alpha_\nu-a_1a_4z_\mu}{a_2a_3}\right)+\beta\right)}}{a_2a_3}\;dz_\mu\eeas
which implies that
\beas f(z)&=&\frac{2\left(a_3K_3-a_1K_4\right)e^{\frac{1}{2}\sum_{j=1}^n \beta_jz_j+\frac{1}{2}\beta}}{\beta_\mu a_2a_3-\beta_\nu a_1a_4}+\alpha_\mu.\eeas
Note that after integration with respect to $z_1$, replacing $\alpha_\nu$ by $a_1a_4z_\mu+a_2a_3z_\nu$, $\alpha_j$ by $z_j$ for $1\leq j(\not=\mu,\nu)\leq n$, where 
$\alpha_j\in\mathbb{C}$ ($1\leq j\leq n$). Hence the solution is 
$\Psi(\alpha_1,\alpha_2,\ldots,\alpha_n)=0$. For simplicity, we suppose 
\bea\label{mn6} f(z)=\frac{2\left(a_3K_3-a_1K_4\right)}{\left(a_2a_3\beta_\mu-a_1a_4\beta_\nu\right)}e^{\frac{1}{2}\sum_{j=1}^n \beta_jz_j+\frac{1}{2}\beta}+h_1\left(y\right),\eea 
where $a_2a_3\beta_\mu-a_1a_4\beta_\nu\not=0$, $a_3K_3-a_1K_4\not=0$, $h_1(y)$ is a finite order entire function in $z_1,z_2,\ldots,z_{\mu-1},a_1a_4z_\mu+a_2a_3z_\nu,z_{\mu+1},\ldots,z_{\nu-1},z_{\nu+1},\ldots,z_n$ with 
$a_2a_3\frac{\pa }{\pa z_\mu}h_1(y)\equiv a_1a_4\frac{\pa }{\pa z_\nu}h_1(y)$.
From (\ref{mn}), with the help of (\ref{poo}) and (\ref{mn6}), we get
\beas h_1(y+s)-h_1(y)\equiv \frac{a_4\beta_\nu K_3-a_2\beta_\mu K_4}{a_2a_3\beta_\mu-a_1a_4\beta_\nu}e^{\frac{1}{2}\sum_{j=1}^n\beta_jz_j+\frac{1}{2}\beta}-\frac{a_2}{a_1}\frac{\pa h_1(y)}{\pa z_\mu}.\eeas
{\bf Sub-case 1.1.2.} Let $g_1(z)$ be independent of both $z_\mu$ and $z_\nu$, then, we have $g(z)=\sum_{j=1}^n \beta_jz_j+g_1(z)+\beta$, where $\beta_j,\beta\in\mathbb{C}$ $(1\leq j\leq n)$ and $g_1(z)$ is a polynomial defined in (\ref{K1}). Using similar arguments to those in \textrm{Sub-case 1.1.1}, we have 
\beas f(z)=\frac{2\left(a_3K_3-a_1K_4\right)}{\left(a_2a_3\beta_\mu-a_1a_4\beta_\nu\right)}e^{\frac{1}{2}\sum_{j=1}^n \beta_jz_j+\frac{1}{2}g_1(z)+\frac{1}{2}\beta}+h_2\left(y\right),\eeas
where $a_2a_3\beta_\mu-a_1a_4\beta_\nu\not=0$, $a_3K_3-a_1K_4\not=0$ with (\ref{K2}) and $h_2(y)$ is a finite order entire function in $z_1,z_2,\ldots,z_{\mu-1},a_1a_4z_\mu+a_2a_3z_\nu,z_{\mu+1},\ldots,z_{\nu-1},z_{\nu+1},\ldots,z_n$ satisfying
$a_2a_3\frac{\pa }{\pa z_\mu} h_2(y)\equiv a_1a_4\frac{\pa}{\pa z_\nu} h_2(y)$ and 
\beas h_2(y+s)-h_2(y)\equiv \frac{a_4\beta_\nu K_3-a_2\beta_\mu K_4}{a_2a_3\beta_\mu-a_1a_4\beta_\nu}e^{\frac{1}{2}\sum_{j=1}^n\beta_jz_j+\frac{1}{2}g_1(z)+\frac{1}{2}\beta}-\frac{a_2}{a_1}\frac{\pa h_2(y)}{\pa z_\mu}.\eeas
{\bf Sub-case 1.2.} When $a_3K_3-a_1K_4=0$, then we have $K_3=a_1/\sqrt{a_1^2+a_3^2}$ and $K_4=a_3/\sqrt{a_1^2+a_3^2}$. Using arguments similar to those in \textrm{Sub-case 1.1.}, we have from (\ref{mn0}) that
\bea\label{t1} f(z)=h_3(y)\quad\text{with}\quad a_2a_3\frac{\pa h_3(y)}{\pa z_\mu}\equiv a_1a_4\frac{\pa h_3(y)}{\pa z_\nu},\eea
where $h_2(y)$ is a finite order transcendental entire function. From (\ref{mn}) and (\ref{t1}), we deduce that
\beas a_1(h_3(y+s)-h_3(y))+a_2\frac{\pa h_3(y)}{\pa z_\mu}=\frac{a_1}{\sqrt{a_1^2+a_3^2}}e^{g(z)/2}.\eeas
{\bf Case 2.} Let $e^{\gamma_2(z+c)-\gamma_1(z+c)}$ be non-constant. 
As demonstrated in \textrm{Case 2} of \textrm{Theorem \ref{TT1}}, the same arguments lead to the conclusion that both $\Gamma_2(z)$ and $\Omega_2(z)$ are not identically zero. Since $e^{\gamma_2(z+c)-\gamma_1(z+c)}$ is non-constant and it is easy to see that 
\beas && N\left(r, \Gamma_2(z)e^{\gamma_1(z)-\gamma_1(z+c)}\right)= N\left(r, \Omega_2(z)e^{\gamma_2(z)-\gamma_1(z+c)}\right)\\[2mm]
&&=N\left(r,\frac{\left(a_1+i a_3\right)K_2}{\left(a_1-ia_3\right)K_1}e^{\gamma_2(z+c)-\gamma_1(z+c)}\right)=N\left(r, 0; \Gamma_2(z)e^{\gamma_1(z)-\gamma_1(z+c)}\right)\\[2mm]
&&=N\left(r, 0; \Omega_2(z)e^{\gamma_2(z)-\gamma_1(z+c)}\right)=N\left(r,0; \frac{\left(a_1+i a_3\right)K_2}{\left(a_1-i a_3\right)K_1}e^{\gamma_2(z+c)-\gamma_1(z+c)}\right)\\[2mm]&&
=S\left(r, \frac{\left(a_1+i a_3\right)K_2}{\left(a_1-i a_3\right)K_1}e^{\gamma_2(z+c)-\gamma_1(z+c)}\right).\eeas
In the light of \textrm{Lemma \ref{lem1}}, it follows from (\ref{nb}) that either 
\beas\text{either}\quad \Gamma_2(z)e^{\gamma_1(z)-\gamma_1(z+c)}\equiv 1\quad\text{or}\quad\Omega_2(z) e^{\gamma_2(z)-\gamma_1(z+c)}\equiv 1.\eeas 
Now the following cases arise.\\[2mm]
{\bf Sub-case 2.1.} Let $\Gamma_1(z)e^{\gamma_1(z)-\gamma_1(z+c)}\equiv 1$. From (\ref{nb}), we have
\beas-\left(ia_4\frac{\pa\gamma_2(z)}{\pa z_\nu}+a_2\frac{\pa\gamma_2(z)}{\pa z_\mu}\right)\equiv (a_1+ ia_3)\left(e^{\gamma_2(z+c)-\gamma_2(z)}-1\right).\eeas
Using arguments similar to those presented in \textrm{Sub-case 2.1} of \textrm{Theorem \ref{TT1}}, we deduce that
$\gamma_1(z)=\sum_{j=1}^nb_jz_j+\xi_1(z)+A$ and $\gamma_2(z)=\sum_{j=1}^nd_jz_j+\xi_2(z)+B$, where $b_j,d_j,A,B\in\mathbb{C}$ $(1\leq j\leq n)$ and  
$\xi_k(z)$ $(k=1,2)$ is a polynomial defined in (\ref{K1}). Thus, we have 
\bea\label{h5}&& (a_1-i a_3)\left(e^{\sum_{j=1}^nb_jc_j}-1\right)\equiv 
a_4b_\nu i-a_2b_\mu+\left(a_4 i\frac{\pa \xi_1(z)}{\pa z_\nu}-a_2 \frac{\pa \xi_1(z)}{\pa z_\mu}\right)\\ \text{and}
&&\label{h6} (a_1+i a_3)\left(1-e^{\sum_{j=1}^nd_jc_j}\right)\equiv 
a_4d_\nu i+a_2d_\mu+\left(a_4i\frac{\pa \xi_2(z)}{\pa z_\nu}+a_2\frac{\pa \xi_2(z)}{\pa z_\mu}\right).\eea
{\bf Sub-case 2.1.1.} Let $\xi_k(z)$ contains the variables $z_\mu$ or $z_\nu$ or both the variables $z_\mu$ and $z_\nu$ with $a_4 i\frac{\pa }{\pa z_\nu}\xi_k(z)+(-1)^k a_2 \frac{\pa }{\pa z_\mu}\xi_k(z)\not=0$, then by comparing the degrees on both sides of (\ref{h5}) and (\ref{h6}), we get that $\deg(\xi_k(z))\leq 1$ for $k=1,2$. For simplicity, 
we still denote $\gamma_1(z)=\sum_{j=1}^nb_jz_j+b_{n+1}$ and $\gamma_2(z)=\sum_{j=1}^nd_jz_j+d_{n+1}$, where $b_j,d_j\in\mathbb{C}$ $(1\leq j\leq n+1)$. 
Therefore, we have
\bea\label{nb1}e^{\sum_{j=1}^nb_jc_j}-1\equiv \frac{a_2b_\mu-a_4b_\nu i}{a_3i-a_1}\quad \text{and}\quad e^{\sum_{j=1}^nd_jc_j}-1\equiv -\frac{a_4d_\nu i+a_2d_\mu}{a_1+a_3i}.\eea
From (\ref{po3}), we have 
\be\label{nb3}a_2a_3\frac{\pa f(z)}{\pa z_\mu}-a_1a_4\frac{\pa f(z)}{\pa z_\nu}=\left(\frac{a_3}{2}-\frac{a_1}{2i}\right)K_1e^{\sum_{j=1}^nb_jz_j+A}+\left(\frac{a_3}{2}+\frac{a_1}{2i}\right)K_2e^{\sum_{j=1}^nd_jz_j+B}.\ee
From (\ref{po1}), (\ref{nb1}) and (\ref{nb3}), using arguments similar to those presented in \textrm{Sub-case 1.1.1} of \textrm{Theorem \ref{TT1}}, we deduce the following
\bea\label{qa} f(z)=\frac{\left(a_3i-a_1\right)K_1e^{\sum_{j=1}^n b_jz_j+A}}{2i\left(a_2a_3b_\mu-a_1a_4b_\nu \right)}+\frac{\left(a_3i+a_1\right)K_2e^{\sum_{j=1}^n d_jz_j+B}}{2i\left(a_2a_3d_\mu-a_1a_4d_\nu\right)}+h_4(y),\eea
where $a_2a_3b_\mu-a_1a_4b_\nu\not=0$, $a_2a_3d_\mu-a_1a_4d_\nu\not=0$, $h_4(y)$ is a finite order entire function with 
$a_2a_3\frac{\pa}{\pa z_\mu}h_4(y)\equiv a_1a_4\frac{\pa }{\pa z_\nu}h_4(y)$. 
Using (\ref{nb1}) and (\ref{qa}), we get from (\ref{po1}) that 
\beas h_4(y+s)-h_4(y)\equiv -\frac{a_2}{a_1}\frac{\pa h_4(y)}{\pa z_\mu}\equiv -\frac{a_4}{a_3}\frac{\pa h_4(y)}{\pa z_\nu}.\eeas
{\bf Sub-case 2.1.2.} Let $\xi_k(z)$ be independent of both $z_\mu$ and $z_\nu$, then, we have $\gamma_1(z)=\sum_{j=1}^nb_jz_j+\xi_1(z)+A$ and $\gamma_2(z)=\sum_{j=1}^nd_jz_j+\xi_2(z)+B$, where $b_j, d_j, A, B\in\mathbb{C}$ $(1\leq j\leq n)$ and  
$\xi_k(z)$ $(k=1,2)$ is a polynomial defined in (\ref{K1}). Using similar arguments to those in \textrm{Sub-case 2.1.1}, we have 
\bea\label{qa} f(z)=\frac{\left(a_3i-a_1\right)K_1e^{\sum_{j=1}^n b_jz_j+\xi_1(z)+A}}{2i\left(a_2a_3b_\mu-a_1a_4b_\nu \right)}+\frac{\left(a_3i+a_1\right)K_2e^{\sum_{j=1}^n d_jz_j+\xi_2(z)+B}}{2i\left(a_2a_3d_\mu-a_1a_4d_\nu\right)}+h_5(y),\eea
where $a_2a_3b_\mu-a_1a_4b_\nu\not=0$, $a_2a_3d_\mu-a_1a_4d_\nu\not=0$ with (\ref{nb1}) and  $h_5(y)$ is a finite order entire function satisfying
$a_2a_3\frac{\pa }{\pa z_\mu}h_5(y)\equiv a_1a_4\frac{\pa}{\pa z_\nu}h_5(y)$ and
\beas h_5(y+s)-h_5(y)\equiv -\frac{a_2}{a_1}\frac{\pa h_5(y)}{\pa z_\mu}\equiv -\frac{a_4}{a_3}\frac{\pa h_5(y)}{\pa z_\nu}.\eeas
{\bf Sub-case 2.2.} Let $\Omega_2(z)e^{\gamma_2(z)-\gamma_1(z+c)}\equiv 1$. From (\ref{nb}), we have
\beas -\frac{K_1}{K_2}\frac{i a_4\frac{\pa\gamma_1(z)}{\pa z_2}-a_2\frac{\pa\gamma_1(z)}{\pa z_1}+\left(a_1-i a_3\right)}{a_1+i a_3}e^{\gamma_1(z)-\gamma_2(z+c)}\equiv 1.\eeas
Applying the same reasoning as in \textrm{Sub-case 2.2} of \textrm{Theorem \ref{TT1}} leads to a contradiction.
This completes the proof.\end{proof} 
\begin{proof}[\bf Proof of the Theorem \ref{TT3}] Let $f$ be a finite order transcendental entire function on $\mathbb{C}^n$ that satisfies (\ref{e3}), where $g(z)$ is a non-constant polynomial on $\mathbb{C}^n$. Using arguments similar to those in \textrm{Theorem \ref{TT1}}, we obtain the following
\be\label{ds1}a_1f(z+c)=\frac{K_1e^{\gamma_1(z)}+K_2e^{\gamma_2(z)}}{2}\;\text{and}\;
a_2\frac{\pa f(z)}{\pa z_\mu}+a_3\frac{\pa^2 f(z)}{\pa z_\mu^2}=\frac{K_1e^{\gamma_1(z)}-K_2e^{\gamma_2(z)}}{2i},\ee
where $K_1,K_2\in\mathbb{C}\setminus\{0\}$ such that $K_1K_2=1$, $\gamma_1(z)=P(z)+g(z)/2$, $\gamma_2(z)=-P(z)+g(z)/2$ and $P(z)$ is a polynomial on $\mathbb{C}^n$.
From (\ref{ds1}), we have
\bea\label{ds3} &&\Gamma_3(z)e^{\gamma_1(z)-\gamma_1(z+c)}+\Omega_3(z)e^{\gamma_2(z)-\gamma_1(z+c)}+\frac{K_2}{K_1}e^{\gamma_2(z+c)-\gamma_1(z+c)}\equiv 1,\eea
where \beas\begin{array}{lll}
\Gamma_3(z)=i\left(\frac{a_3}{a_1}\left(\frac{\pa^2\gamma_1(z)}{\pa z_\mu^2}+\left(\frac{\pa\gamma_1(z)}{\pa z_\mu}\right)^2\right)+\frac{a_2}{a_1}\frac{\pa\gamma_1(z)}{\pa z_\mu}\right),\\
\\
\Omega_3(z)=\frac{iK_2}{K_1}\left(\frac{a_3}{a_1}\left(\frac{\pa^2\gamma_2(z)}{\pa z_\mu^2}+\left(\frac{\pa\gamma_2(z)}{\pa z_\mu}\right)^2\right)+\frac{a_2}{a_1}\frac{\pa\gamma_2(z)}{\pa z_\mu}\right).\end{array}\eeas
It is necessary to consider the following cases separately.\\[2mm]
{\bf Case 1.} Let $e^{\gamma_2(z+c)-\gamma_1(z+c)}$ be constant. By means of arguments similar to those in \textrm{Case 1} of \textrm{Theorem \ref{TT1}}, we deduce that
$P(z+c)\equiv-k/2$, where $k\in\mathbb{C}$.
From (\ref{ds1}), we have
\bea\label{ds4}a_1f(z+c)=K_3e^{\frac{g(z)}{2}}\;\;\text{and}\;\;a_2\frac{\pa f(z)}{\pa z_\mu}+a_3\frac{\pa^2 f(z)}{\pa z_\mu^2}=K_4e^{\frac{g(z)}{2}},\eea
where $K_3=\frac{K_1\rho+K_2\rho^{-1}}{2}$, $K_4=\frac{K_1\rho-K_2\rho^{-1}}{2i}$, $e^{p}=\rho\in\mathbb{C}\setminus\{0\}$ and $K_3^2+K_4^2=1$. Clearly $K_3\not=0$. Now the following cases arise separately.\\
{\bf Sub-case 1.1.} When $K_4=0$. Then from (\ref{ds4}), we deduce that $f(z)=\frac{K_3}{a_1}e^{\frac{g(z-c)}{2}}$, where $g(z)$ is a non-constant polynomial on $\mathbb{C}^n$ with $2a_2\frac{\pa }{\pa z_\mu}g(z)+a_3\left(\frac{\pa}{\pa z_\mu}g(z)\right)^2+2a_3\frac{\pa^2}{\pa z_\mu^2}g(z)\equiv 0$.\\
{\bf Sub-case 1.2.} When $K_4\not=0$.
From (\ref{ds4}), we deduce that 
\bea\label{ds5}\frac{K_3}{K_4}\left(\frac{a_2}{2a_1}\frac{\pa g(z)}{\pa z_\mu}+\frac{a_3}{4a_1}\left(\left(\frac{\pa g(z)}{\pa z_\mu}\right)^2+2\frac{\pa^2g(z)}{\pa z_\mu^2}\right)\right)\equiv e^{\frac{g(z+c)-g(z)}{2}}.\eea
From (\ref{ds5}), using arguments similar to those presented in \textrm{Case 1} of \textrm{Theorem \ref{TT1}}, we deduce that
$g(z)=\sum_{j=1}^n \beta_jz_j+g_1(z)+\beta$, where $\beta_j, \beta\in\mathbb{C}$ $(1\leq j\leq n)$ and $g_1(z)$ is a polynomial defined in (\ref{K1}). From (\ref{ds5}), we have
\bea\label{j1}\frac{K_3}{K_4}\left(\frac{a_2}{2a_1}\left(\beta_\mu+\frac{\pa g_1(z)}{\pa z_\mu}\right)+\frac{a_3}{4a_1}\left(\left(\beta_\mu+\frac{\pa g_1(z)}{\pa z_\mu}\right)^2+2\frac{\pa^2 g_1(z)}{\pa z_\mu^2}\right)\right)\equiv e^{\frac{1}{2}\sum_{j=1}^n \beta_jc_j}.\eea
{\bf Sub-case 1.2.1.} When $g_1(z)$ contain the variable $z_\mu$, then by comparing the degrees on both sides of (\ref{j1}), we get that $\deg(g_1(z))\leq 1$. For simplicity, we still denote 
$g(z)=\sum_{j=1}^n\beta_jz_j+\beta$, where $\beta_j,\beta\in\mathbb{C}$ $(1\leq j\leq n)$. This implies that $g_1(z)\equiv 0$. From (\ref{ds4}) and (\ref{ds5}), we have  
\beas f(z)=\frac{K_3}{a_1}e^{\frac{1}{2}\sum_{j=1}^n \beta_j(z_j-c_j)+\frac{1}{2}\beta}\quad\text{and}\quad e^{\frac{1}{2}\sum_{j=1}^n \beta_jc_j}\equiv \frac{K_3}{K_4}\left(\frac{a_2}{2a_1}\beta_\mu+\frac{a_3}{4a_1}\beta_\mu^2\right),\eeas
where $\beta_j,\beta,K_3(\not=0),K_4(\not=0)\in\mathbb{C}$ ($1\leq j\leq n$) with $\beta_\mu\not=0$ and $K_3^2+K_4^2=1$.\\[2mm]
{\bf Sub-case 1.2.2.} If $g_1(z)$ is independent of $z_\mu$, then we have $g(z)=\sum_{j=1}^n\beta_jz_j+g_1(z)+\beta$, where $\beta_j,\beta\in\mathbb{C}$ $(1\leq j\leq n)$. From (\ref{ds4}) and (\ref{ds5}), we have  
\beas f(z)=\frac{K_3}{a_1}e^{\frac{1}{2}\sum_{j=1}^n \beta_j(z_j-c_j)+\frac{1}{2}g_1\left(z\right)+\frac{1}{2}\beta}, e^{\frac{1}{2}\sum_{j=1}^n \beta_jc_j}\equiv \frac{K_3}{K_4}\left(\frac{a_2}{2a_1}\beta_\mu+\frac{a_3}{4a_1}\beta_\mu^2\right).\eeas
{\bf Case 2.} Let $e^{\gamma_2(z+c)-\gamma_1(z+c)}$ be non-constant. As demonstrated in \textrm{Case 2} of \textrm{Theorem \ref{TT1}}, the same arguments lead to the conclusion that $\Gamma_3(z)\not\equiv 0$ and $\Omega_3(z)\not\equiv 0$.
As $e^{\gamma_2(z+c)-\gamma_1(z+c)}$ is non-constant, it is evident that 
\beas&& N\left(r, \Gamma_3(z)e^{\gamma_1(z)-\gamma_1(z+c)}\right)= N\left(r, \Omega_3(z)e^{\gamma_2(z)-\gamma_1(z+c)}\right)=N\left(r,K_2e^{\gamma_2(z+c)-\gamma_1(z+c)}/K_1\right)\\[2mm]
&&=N\left(r, 0; \Gamma_3(z)e^{\gamma_1(z)-\gamma_1(z+c)}\right)=N\left(r, 0; \Omega_3(z)e^{\gamma_2(z)-\gamma_1(z+c)}\right)\\[2mm]
&&=N\left(r,0; K_2e^{\gamma_2(z+c)-\gamma_1(z+c)}/K_1\right)=S\left(r, K_2e^{\gamma_2(z+c)-\gamma_1(z+c)}/K_1\right).\eeas
In the light of \textrm{Lemma \ref{lem1}}, it follows from (\ref{ds3}) that either 
\beas \text{either} \quad\Gamma_3(z)e^{\gamma_1(z)-\gamma_1(z+c)}\equiv 1\quad \text{or} \quad\Omega_3(z) e^{\gamma_2(z)-\gamma_1(z+c)}\equiv 1.\eeas 
Now the following cases arise.\\ [2mm]
{\bf Sub-case 2.1.} Let $\Gamma_3(z)e^{\gamma_1(z)-\gamma_1(z+c)}\equiv 1$. From (\ref{ds3}), we have 
\bs\beas -\frac{K_1}{K_2}\Omega_3(z)e^{\gamma_2(z)-\gamma_2(z+c)}\equiv 1. \eeas\es
Using arguments similar to those presented in \textrm{Sub-case 2.1} of \textrm{Theorem \ref{TT1}}, we deduce that
$\gamma_1(z)=\sum_{j=1}^nb_jz_j+\xi_1(z)+A$ and $\gamma_2(z)=\sum_{j=1}^nd_jz_j+\xi_2(z)+B$, where $b_j,d_j,A,B\in\mathbb{C}$ ($1\leq j\leq n$) and  
$\xi_k(z)$ $(K=1,2)$ is a polynomial defined in (\ref{K1}). Therefore, we have 
\bea\label{j2}\begin{array}{lll}
-ia_1e^{\sum_{j=1}^nb_jc_j}\equiv a_3\left(\frac{\pa^2 \xi_1(z)}{\pa z_\mu^2}+\left(b_\mu+\frac{\pa \xi_1(z)}{\pa z_\mu})^2\right)+a_2\left(b_\mu+\frac{\pa \xi_1(z)}{\pa z_\mu}\right)\right),\\\\
ia_1e^{\sum_{j=1}^nd_jc_j}\equiv a_3\left(\frac{\pa^2 \xi_2(z)}{\pa z_\mu^2}+\left(d_\mu+\frac{\pa \xi_2(z)}{\pa z_\mu})^2\right)+a_2\left(d_\mu+\frac{\pa \xi_2(z)}{\pa z_\mu}\right)\right).\end{array}\eea
{\bf Sub-case 2.1.1.} If $\xi_k(z)$ is dependent on the variable $z_\mu$, then by comparing the degrees on both sides of (\ref{j2}), we get that $\deg(\xi_k(z))\leq 1$ for $k=1,2$. For simplicity, 
we still denote $\gamma_1(z)=\sum_{j=1}^nb_jz_j+A$ and $\gamma_2(z)=\sum_{j=1}^nd_jz_j+B$, where $b_j,d_j,A,B\in\mathbb{C}$ $(1\leq j\leq n)$. 
Thus, we conclude that $\xi_k(z)\equiv 0$ for $k=1,2$.
Therefore, we have 
\bea\label{j3} e^{\sum_{j=1}^nb_jc_j}\equiv i\left(a_3b_\mu^2+a_2b_\mu\right)/a_1\quad\text{and}\quad e^{\sum_{j=1}^nd_jc_j}\equiv -i\left(a_3d_\mu^2+a_2d_\mu\right)/a_1,\eea
which implies that $b_\mu\not=0$, $d_\mu\not=0$. From (\ref{ds1}), we derive that
\beas f(z)=\frac{K_1e^{\sum_{j=1}^nb_j(z_j-c_j)+A}+K_2e^{\sum_{j=1}^nd_j(z_j-c_j)+B}}{2a_1},\eeas 
where $b_j,d_j,A,B,K_1,K_2\in\mathbb{C}$ ($1\leq j\leq n$) with $b_\mu\not=0$, $d_\mu\not=0$ and $K_1K_2=1$.\\[2mm]
{\bf Sub-case 2.1.2.} If $\xi_k(z)$ $(k=1,2)$ is independent of $z_\mu$, then we have 
$\gamma_1(z)=\sum_{j=1}^nb_jz_j+\xi_1\left(z\right)+A$ and $\gamma_2(z)=\sum_{j=1}^nd_jz_j+\xi_2\left(z\right)+B$,
where $b_j,d_j,A,B\in\mathbb{C}$ $(1\leq j\leq n)$ and $\xi_k(z)$ $(k=1,2)$ is a polynomial defined in (\ref{K1}).
From (\ref{ds1}) and (\ref{j2}), we deduce (\ref{j3}) and 
\beas f(z)=\frac{K_1e^{\sum_{j=1}^nb_j(z_j-c_j)+\xi_1\left(z\right)+A}+K_2e^{\sum_{j=1}^nd_j(z_j-c_j)+\xi_2\left(z\right)+B}}{2a_1},\eeas 
where $b_j,d_j,A,B,K_1,K_2\in\mathbb{C}$ ($1\leq j\leq n$) with $b_\mu\not=0$, $d_\mu\not=0$ and $K_1K_2=1$.\\[2mm]
{\bf Sub-case 2.2.} Let $\Omega_3(z)e^{\gamma_2(z)-\gamma_1(z+c)}\equiv 1$. From (\ref{ds3}), we get 
$-K_1\Gamma_3(z)e^{\gamma_1(z)-\gamma_2(z+c)}/K_2\equiv 1$. Applying the same reasoning as in \textrm{Sub-case 2.2} of \textrm{Theorem \ref{TT1}} leads to a contradiction. This completes the proof.
\end{proof} 
\begin{proof}[\bf Proof of the Theorem \ref{TT4}] Let $f$ be a finite order transcendental entire function on $\mathbb{C}^n$ satisfies (\ref{e4}). Using arguments similar to those in \textrm{Theorem \ref{TT1}}, we obtain the following
\be\label{dds1}a_1f(z+c)=\frac{K_1e^{\gamma_1(z)}+K_2e^{\gamma_2(z)}}{2}\quad\text{and}\quad
a_2\frac{\pa f(z)}{\pa z_\mu}+a_3\frac{\pa^2 f(z)}{\pa z_\mu \pa z_\nu}=\frac{K_1e^{\gamma_1(z)}-K_2e^{\gamma_2(z)}}{2i},\ee
where $K_1,K_2\in\mathbb{C}\setminus\{0\}$ such that $K_1K_2=1$, $\gamma_1(z)=P(z)+g(z)/2$, $\gamma_2(z)=-P(z)+g(z)/2$ and $P(z)$ is a polynomial on $\mathbb{C}^n$. From (\ref{dds1}), we deduce
\bea\label{dds3}\Gamma_4(z)e^{\gamma_1(z)-\gamma_1(z+c)}+\Omega_4(z)e^{\gamma_2(z)-\gamma_1(z+c)}+\frac{K_2}{K_1}e^{\gamma_2(z+c)-\gamma_1(z+c)}\equiv 1,\eea
where \beas\begin{array}{lll}
\Gamma_4(z)=i\left(\frac{a_3}{a_1}\left(\frac{\pa^2\gamma_1(z)}{\pa z_\mu\pa z_\nu}+\frac{\pa\gamma_1(z)}{\pa z_\mu}\frac{\pa\gamma_1(z)}{\pa z_\nu}\right)+\frac{a_2}{a_1}
\frac{\pa\gamma_1(z)}{\pa z_\mu}\right),\\\\
\Omega_4(z)=\frac{iK_2}{K_1}\left(\frac{a_3}{a_1}\left(\frac{\pa^2\gamma_2(z)}{\pa z_\mu\pa z_\nu}+\frac{\pa\gamma_2(z)}{\pa z_\mu}
\frac{\pa\gamma_2(z)}{\pa z_\nu}\right)+\frac{a_2}{a_1}\frac{\pa\gamma_2(z)}{\pa z_\mu}\right).\end{array}\eeas
The rest of the proof follows using the same arguments as in the proof of \textrm{Theorem \ref{TT3}}. Thus, the conclusions of this theorem are straightforward. This completes the proof.
\end{proof} 
\section{Declarations}
\noindent{\bf Acknowledgments:} The third author is supported by a grant from the University Grants Commission (IN) (No. F. 44-1/2018 (SA-III)). We would also want to thank the anonymous reviewers and the editing team for their suggestions.\\
{\bf Conflict of Interest:} The authors declare that we do not have any conflicts of interest.\\
{\bf Data availability:} Not applicable.

\end{document}